\documentclass[11pt,oneside]{article}

\usepackage[round,colon,authoryear]{natbib} 

\usepackage{graphicx}
\usepackage{amsmath}
\usepackage{amsfonts}
\usepackage{amssymb}
\usepackage{amsthm}
\usepackage{enumerate}
\usepackage{pdfsync}
\usepackage{times}

     \usepackage[final,notref,notcite,color]{showkeys}

\usepackage{bm}

    \numberwithin{equation}{section}

\newcommand{\R}{\mathbb{R}}            
\newcommand{\N}{\mathbb{N}}            

\newcommand{\al}{\alpha}                
\newcommand{\ga}{\gamma}                
\newcommand{\dx}{\mathrm{d}}    
\newcommand{\ep}{\varepsilon}              

\newcommand{\ra}{\rightarrow}           
\newcommand{\da}{\downarrow}           

\newcommand{\Om}{\Omega}

\bmdefine\mub{\mu}
\bmdefine\etab{\eta}
\bmdefine\thetab{\vartheta}
\bmdefine\betab{\beta}
\bmdefine\sigmab{\sigma}
\bmdefine\varsigmab{\varsigma}
\bmdefine\taub{\tau}
\bmdefine\rhob{\rho}
\bmdefine\varrhob{\varrho}
\bmdefine\gammab{\gamma}
\bmdefine\Gammab{\Gamma}

\bmdefine\xib{\mathbf{\xi}}
\bmdefine\etab{\mathbf{\eta}}

\theoremstyle{plain}
\newtheorem{thm}{Theorem}[section]

\newtheorem{prop}[thm]{Proposition}

\theoremstyle{definition}

\newtheorem{exmp}{Example}[section]

\theoremstyle{remark}
\newtheorem{rem}{Remark}[section]

\title{   
\textsc{Skew-Unfolding   the   Skorokhod 
Reflection of a Continuous Semimartingale} \thanks{~ We are indebted to Marcel Nutz and Daniel Ocone for stimulating discussions,  to Vilmos Prokaj and Johannes Ruf for their careful reading of the manuscript and for  their suggestions, and to the referee for simplifying the last part of our argument in Example \ref{Counter}.} 
}

\vspace{2mm}

\author{  
 TOMOYUKI ICHIBA      \thanks{ ~
Department of Statistics and Applied Probability, South Hall, University of California, Santa Barbara, CA 93106, USA (E-mail:    {\it ichiba@pstat.ucsb.edu}). Research supported in part by the National Science Foundation under grant NSF-DMS-13-13373. 
          }
 \and
 IOANNIS KARATZAS \thanks{~
Department of Mathematics,  Columbia University (E-mail: {\it ik@math.columbia.edu}), and       \textsc{Intech} Investment Management,  One Palmer Square, Suite 441, Princeton, NJ 08542, USA    (E-mail:    {\it ik@enhanced.com}). Research   supported in part by  the National Science Foundation under  grants NSF-DMS-09-05754 and NSF-DMS-14-05210.}  
                                      }

\date{July 2, 2014}

\begin{document}

\maketitle

\bigskip

\centerline{
\it Dedicated to  Terry Lyons on the occasion of his  60th  birthday
}

\bigskip

\begin{abstract}

\bigskip \noindent \small  
The  Skorokhod   reflection of a continuous semimartingale is unfolded, in a possibly skewed manner, into another continuous semimartingale on an enlarged probability space according to the excursion-theoretic methodology  of  Prokaj  (2009). This is done in terms of a skew version of the  Tanaka equation, whose properties   are studied in some detail. The result is used to construct a system of two diffusive  particles with rank-based characteristics and skew-elastic collisions. Unfoldings of conventional reflections are also discussed, as are examples involving skew Brownian Motions and skew Bessel processes. 
\end{abstract}

\medskip
\noindent{\it Keywords and Phrases:}  Skorokhod and conventional reflections; skew and perturbed Tanaka equations; skew Brownian and Bessel processes; pure and Ocone martingales; local time; competing particle systems; asymmetric collisions. 

\medskip
\noindent{\it AMS 2000 Subject Classifications:}  Primary, 60G42; secondary, 60H10.

\vspace{1mm}

 \input amssym.def
\input amssym

 \bigskip
\medskip


\section{The Result}
\label{Res}

On a filtered probability space $(\Om, \mathcal{F}, \mathbb{P}), \,\mathbb{F} = \{ \mathcal{F} (t) \}_{0 \le t < \infty}$ satisfying the so-called  ``usual conditions" of right continuity and augmentation by null sets, we consider a  real-valued continuous semimartingale $U(\cdot)$ of the form 
\begin{equation}
\label{1}
U (t) \, =\, M(t) + A(t)\,, \qquad 0 \le t < \infty 
\end{equation}
with $M(\cdot) $ a continuous local martingale  and  $A(\cdot)$ a process of finite first  variation on compact intervals. We   assume $M(0)=A(0)=0$ for concreteness. 

There are two ways to ``fold", or  reflect, this semimartingale about the origin. One is the {\it conventional reflection}
\begin{equation}
\label{2a}
R (t) \, :=\, | U(t) |\,, \qquad 0 \le t < \infty\,;
\end{equation}
the other is the {\it \textsc{Skorokhod} reflection}
\begin{equation}
\label{2}
S(t) \,:=\, U(t) + \max_{0 \le s \le t} \big( - U(s) \big)\,, \qquad 0 \le t < \infty\,.
\end{equation}
The following result, inspired by \textsc{Prokaj} (2009),  shows how  the first   can be obtained from the second, by suitably  unfolding  the \textsc{Skorokhod} reflection in a possibly ``skewed" manner.

\begin{thm}
 \label{Thm1}
Fix a constant $\alpha \in (0,1)$. There exists an enlargement $\,\big(\widetilde{\Om}, \widetilde{\mathcal{F}}, \widetilde{\mathbb{P}}\big),\, \widetilde{\mathbb{F}} = \{ \widetilde{\mathcal{F}} (t) \}_{0 \le t < \infty}\,$ of the filtered probability space $(\Om, \mathcal{F}, \mathbb{P}), \,\mathbb{F} = \{ \mathcal{F} (t) \}_{0 \le t < \infty}$  with a measure-preserving map $\,\pi : \Om \ra \widetilde{\Om}\,$, and on this enlarged space a continuous semimartingale $X(\cdot)$    that satisfies
\begin{equation}
\label{3}
\big|X(\cdot) \big| = S(\cdot)\,,\qquad L^X (\cdot) =\al \,  L^S (\cdot)\,, \qquad X(\cdot) \,=\int_0^{\, \cdot} \overline{\text{sgn}} \big( X(t) \big)\, \dx U(t) + {\, 2 \, \al - 1\, \over \al}\, L^X (\cdot)\,. 
\end{equation}
\end{thm}

\medskip
Here and throughout this note, we  use the notation 
\begin{equation}
\label{LT}
L^U (\cdot) \, :=\, \lim_{\ep \da 0} \, { 1 \over \, 2 \, \ep\,} \int_0^{\, \cdot} 1_{ \{0 \le U(t) < \ep\} }\, \dx  \langle U \rangle (t)\,, \qquad   \widehat{L}^{U}(\cdot) \, :=\,  { 1 \over \,2\,}\, \left(L^{ U}(\cdot) + L^{- U}(\cdot) \right)
\end{equation}
respectively for the {\it right} and the {\it symmetric} local time at the origin of a continuous semimartingale as in (\ref{1}), and the conventions 
$$
\overline{\text{sgn}} (x)\,:=\, \mathbf{ 1}_{(0, \infty)}(x)- \mathbf{ 1}_{(-\infty, 0)} (x)\,, \qquad  \text{sgn} (x) \,:=\, \mathbf{ 1}_{(0, \infty)}(x)- \mathbf{ 1}_{(-\infty, 0]} (x) \,, \qquad x \in \R
$$
for the symmetric and the left-continuous versions, respectively,  of the signum function. We also denote by $\, \mathbb{F}^U = \{ \mathcal{F}^U (t) \}_{0 \le t < \infty}$ the ``natural filtration" of $\, U(\cdot)\,$, that is, the smallest filtration that satisfies the usual conditions  and with respect  to which   $\, U(\cdot)\,$ is adapted; we set $\, \mathcal{F}^U (\infty) := \sigma \big( \bigcup_{\,0 \le t < \infty}  \mathcal{F}^U (t) \big)\,$. Equalities between stochastic processes, such as in (\ref{3}), are to be understood throughout in the almost sure sense. 

Theorem \ref{Thm1} constructs  a continuous semimartingale $X(\cdot)$  whose conventional reflection   coincides with the \textsc{Skorokhod} reflection 
of the given semimartingale $U(\cdot)$, and which satisfies the stochastic integral equation in (\ref{3}). We think of this equation  as a  {\it skew  version} of the celebrated \textsc{Tanaka} equation  driven by the continuous semimartingale $U(\cdot)$,  whose ``skew-unfolding" it produces via the      parameter $ \, \al \, $. When there is no skewness, i.e., with $\al = 1/2\,$, the integral equation of (\ref{3}) reduces to the classical \textsc{Tanaka} equation
; in this case Theorem \ref{Thm1} is just the main result in the paper by \textsc{Prokaj} (2009), which     inspired our work. 

We shall prove Theorem \ref{Thm1} in section \ref{Pf}, then use it in section \ref{Appl}   to construct a system of two diffusive  particles with rank-based characteristics and skew-elastic collisions. Section \ref{CR} discusses a similar skew-unfolding of the conventional reflection $\,R(\cdot) = |U(\cdot)|\,$ of $ \,U(\cdot) \,$. In the section that follows  we discuss briefly some properties of the {\it skew \textsc{Tanaka} equation} in (\ref{3}).

\section{The Skew Tanaka Equation}
\label{ST}

 A first question that arises regarding the stochastic integral equation in (\ref{3}), is whether it can be written in the more conventional form 
 \begin{equation}
\label{3a}
 X(\cdot) \,=\int_0^{\, \cdot}  \text{sgn} \big( X(t) \big)\, \dx U(t) + {\, 2 \, \al - 1\, \over \al}\, L^X (\cdot)\,, 
\end{equation}
in terms of the asymmetric (left-continuous) version of the signum function. 
\newpage
For this, it is necessary and sufficient to have 
\begin{equation}
\label{3b}
\int_0^{\, \cdot}  \mathbf{ 1}_{ \{ X(t) =0\} }\, \dx U(t)\equiv 0\,, \qquad \text{or equivalently} \qquad \int_0^{\, \cdot}  \mathbf{ 1}_{ \{ S(t) =0\} }\, \dx U(t)\equiv 0 \end{equation}
in the context of Theorem \ref{Thm1}. Now from (\ref{1}), (\ref{2}) it is clear that $M(\cdot)$ is the local martingale part of the continuous semimartingale $S(\cdot)$, so we have $\,  \langle S \rangle (\cdot)  = \langle U \rangle (\cdot) = \langle M \rangle (\cdot) \,$ and 
\begin{equation}
\label{3c}
\int_0^{\infty}  \mathbf{ 1}_{ \{ S(t) =0\} }\, \dx \langle M \rangle (t)\,=\, 0 
\end{equation}
(e.g., \textsc{Karatzas \& Shreve}, Exercise 3.7.10). This gives $\,\int_0^{\, \cdot}  \mathbf{ 1}_{ \{ S(t) =0\} }\, \dx M(t)\equiv 0\,$,  so   (\ref{3b}) will follow if and only if   
\begin{equation}
\label{3cc}
\int_0^{\, \cdot} \mathbf{ 1}_{ \{ S(t) =0\} }\, \dx A(t)\,\equiv\, 0 
\end{equation}
holds; and on the strength of (\ref{3c}), a sufficient condition for (\ref{3cc})   is that  $A(\cdot)$ be absolutely continuous with respect to the quadratic variation process  $\, \langle M \rangle (\cdot)$. We have the following result. 


\begin{prop}
 \label{Prop1}
For a given continuous semimartingale $U(\cdot)$ of the form (\ref{1}) the  stochastic integral equation of (\ref{3}) can be cast equivalently in the form (\ref{3a}), if and only if  (\ref{3cc}) holds;  and in this case we have the identification $\, L^S (t) = \max_{\,0 \le s \le t} \big( - U(s) \big)\,$ and the filtration comparisons
\begin{equation}
\label{FiltComp}
{\cal F}^{|X|} (t)= {\cal F}^U (t) \subseteq {\cal F}^{X} (t)\,, \quad 0 \le t < \infty\,.
\end{equation}
 Whereas, a sufficient condition for (\ref{3cc}) to hold, is that   there exist  an $\mathbb{F}-$progressively measurable process $\,p(\cdot)\,$, locally integrable with respect to $\, \langle M \rangle (\cdot)\,$  and such that 
\begin{equation}
\label{3d}
A(\cdot)\,=\int_0^{\, \cdot} p(t)\, \dx \langle M \rangle (t)\,.
\end{equation}
\end{prop}

\noindent
{\it Proof:} The first and third claims have already been argued. As for the second, we observe that   the \textsc{It\^o-Tanaka} formula applied to (\ref{3a}) gives
$$
 S(\cdot)  =  \big| X(\cdot) \big|   =\int_0^{\, \cdot}  \text{sgn} \big( X(t) \big)\, \dx X(t) +   2 \,  L^X (\cdot) =U(\cdot) - {\, 2 \, \al - 1\, \over \al}\, L^X (\cdot)+   2 \,  L^X (\cdot) =  U(\cdot)  +      L^S (\cdot)
$$
on the strength of the second equality in (\ref{3}). It is clear from this expression that the filtration comparison $\, {\cal F}^U (t) \subseteq {\cal F}^{S} (t)\,$ holds for all $\, 0 \le t < \infty\,$; whereas the reverse inclusion and the claimed identification are direct consequences of (\ref{2}). \qed

\smallskip
\noindent
{\it Remark:} More generally (that is, in the  absence of condition (\ref{3cc})), the local time at the origin of the \textsc{Skorokhod} reflection $\,S(\cdot)\,$ is   $\, L^S (t) = \max_{\,0 \le s \le t} \, ( - U(s)  )+ \int_0^{ t} \mathbf{ 1}_{ \{ S(u) =0\} }\, \dx A(u)\,, ~~0 \le t < \infty\,.$  

\subsection{Uniqueness in Distribution for the Skew Tanaka Equation}
\label{Uni}

A second question that arises regarding the skew-\textsc{Tanaka} equation of (\ref{3}), is whether it can be solved uniquely. It is   well-known that we cannot expect pathwise uniqueness or strength to hold for this equation. Such strong existence and uniqueness  fail  already with $\, \al = 1/2\,$ and $\, U(\cdot)\,$ a standard Brownian motion, in which case we have  in (\ref{FiltComp}) also the strict inclusion   $\, {\cal F}^U (t) \subsetneqq {\cal F}^X (t)\,$ for all $\, t \in (0, \infty)\,$ (e.g., \textsc{Karatzas \& Shreve} (1991), Example 5.3.5). The \textsc{Skorokhod} reflection of $\,U(\cdot)\,$ can then be ``unfolded" into a Brownian motion $X(\cdot)$,  whose filtration is strictly finer  than that of the original Brownian motion $\,U(\cdot)$: the unfolding cannot be accomplished without the help of some additional randomness. 
\newpage
 
The issue, therefore, is whether   {\it uniqueness in distribution} holds for the skew-\textsc{Tanaka} equation of (\ref{3}), under appropriate conditions.  We shall address this question in the case of a  continuous local martingale $\,U(\cdot)\,$  with $\, U(0)=0\,$ and  $\,\langle U \rangle (\infty) = \infty\,$. Let us recall a few notions and facts   about such a process, starting with its \textsc{Dambis-Dubins-Schwarz} representation
\begin{equation}
\label{Ocone}
U( t)= B \big( \langle U \rangle ( t) \big)\,, \qquad 0 \le t < \infty
\end{equation}
(cf.$\,$\textsc{Karatzas \& Shreve} (1991), Theorem 3.4.6); here $\, B (\theta) = U (Q (\theta)), \, 0 \le \theta <\infty\,$ is standard Brownian motion, and $\, Q (\cdot)\,$ the right-continuous inverse of the continuous, increasing process $\, \langle U \rangle (\cdot)$. 

We say that this $\,U(\cdot)\,$ is {\it pure,}  if each $\, \langle U \rangle (t)\,$ is   $\, \mathcal{F}^{B}(\infty)-$measurable; we say that it is an  {\it \textsc{Ocone} martingale,} if the processes $\, B(\cdot)\,$ and $\, \langle U \rangle (\cdot)$ are independent (cf.$\,$\textsc{Ocone}  (1993) and \textsc{Dubins et al.}$\,$(1993), Appendix). As discussed   in     \textsc{Vostrikova \& Yor} (2000), a pure   \textsc{Ocone}  martingale  is a Gaussian process.

\begin{prop}
 \label{Prop2}
Suppose that $\, U(\cdot)\,$ is a continuous local martingale with $\, U(0)=0\,$ and  $\,\langle U \rangle (\infty) = \infty\,$. Then uniqueness in distribution holds for the skew-\textsc{Tanaka} equation of (\ref{3}), or equivalently of  (\ref{3a}), provided that either

\smallskip
\noindent
(i) $\, U(\cdot)\,$ is pure; or that 

\smallskip
\noindent
(ii) the quadratic variation process  $\, \langle U \rangle (\cdot)\,$   is adapted to a Brownian motion $\, \Gamma (\cdot) := ( \Gamma_1 (\cdot), \cdots, \Gamma_n (\cdot) )' $, with values in  some Euclidean space $\,  \R^n\,$  and independent of the real-valued  Brownian motion $\, B(\cdot)\,$ in the representation (\ref{Ocone}). 
\end{prop}

\noindent
{\it Proof:} Let us consider a continuous local martingale $  U(\cdot) $ with $U(0)=0$
, and any continuous semimartingale $X(\cdot)$ that satisfies    the stochastic integral equation in (\ref{3}). Then $X(\cdot)$ also satisfies the equation of (\ref{3a}), as the condition (\ref{3d}) holds in this case trivially with $\, p (\cdot) \equiv 0$. In fact, the equation (\ref{3a})    can be written then in the form
$$
X(Q (s)) \,=\int_0^{s}  \text{sgn} \big( X(Q (\theta)) \big)\, \dx B(\theta) + {\, 2 \, \al - 1\, \over \al}\, L^X (Q (s))\,, \qquad 0 \le s < \infty\,,
$$
with $\, Q (\cdot)\,$ the right-continuous inverse of the continuous, increasing process $\, \langle U \rangle (\cdot)$; cf. Proposition 3.4.8 in \textsc{Karatzas \& Shreve} (1991). Setting 
$$
\widetilde{X} (s) := X \big( Q (s) \big)\,, \quad \text{it is straightforward to check} \quad L^{\widetilde{X}} (s) = L^X \big( Q (s)\big)\,, \quad 0 \le s < \infty\,;
$$
for this, one uses the representation (\ref{LT}) for the local time at the origin, along with the fact that the local martingale part of the continuous seminartingale $\, X(\cdot)\,$ in (\ref{3a}) has quadratic variation process $\, \langle U \rangle (\cdot)\,$. Thus, the time-changed process $\, \widetilde{X} (\cdot)\,$ satisfies the stochastic integral equation 
\begin{equation} 
\label{HS0}
\widetilde{X}(s) \,=\int_0^{s}  \text{sgn} \big( \widetilde{X} (\theta)  \big)\, \dx B(\theta) + {\, 2 \, \al - 1\, \over \al}\, L^{\widetilde{X}} (s)\,, \qquad 0 \le s < \infty\,.
\end{equation}
This can be cast as the \textsc{Harrison-Shepp} (1981) equation 
\begin{equation}
\label{HS}
\widetilde{X}(\cdot) \,=\,   \widetilde{W}(\cdot) + {\, 2 \, \al - 1\, \over \al}\, L^{\widetilde{X}} (\cdot) 
\end{equation}
for the skew Brownian motion,  driven by the standard Brownian motion 
\begin{equation}
\label{tildeW}
  \widetilde{W}(\cdot) := \int_0^{\, \cdot}  \text{sgn} \big( \widetilde{X} (\theta)  \big)\, \dx B(\theta)\,.
\end{equation}
It is well-known from the theory of \textsc{Harrison \& Shepp} (1981) that the equation (\ref{HS}) has a pathwise unique, strong solution; in fact, the   skew Brownian motion  
 $ \widetilde{X}(\cdot)$ and the  Brownian motion $ \,\widetilde{W}(\cdot)$ generate   the same filtration. Since 
 \begin{equation}
\label{Rep}
X(t) \,=\, \widetilde{X}   \big( \langle U \rangle ( t) \big)\,, \qquad 0 \le t < \infty  
\end{equation}
holds with    $\, \widetilde{X} (\cdot)\,$  adapted to $\, \mathbb{F}^{\,\widetilde{W}}$, the distribution of $\, X(\cdot)\,$ is uniquely determined whenever 
\begin{equation}
\label{claim}
\text{the Brownian motion}~ \widetilde{W}(\cdot)~ \text{of (\ref{tildeW}) is independent of the   process} ~\langle U \rangle (\cdot)\,, 
\end{equation}
or whenever
\begin{equation}
\label{"pure"}
\langle U \rangle (t) ~\,~\text{is}~\,\, \mathcal{F}^{\,\widetilde{W}}(\infty)-\text{measurable, for every} ~  ~t \in [0, \infty)\,.
\end{equation}

\medskip
But (\ref{"pure"})    holds when $\, U(\cdot)\,$ is  pure (case {\it (i)} of the Proposition); this is   because  from (\ref{tildeW}) we have   $\, B (\cdot)  = \int_0^{\, \cdot}  \text{sgn} \big( \widetilde{X} (\theta)  \big)\, \dx \widetilde{W}(\theta)\,$,     therefore   $ \, {\cal F}^B (t) \subseteq {\cal F}^{\widetilde{W}} (t) $ for all $\, t \in [0, \infty)$ and thus   $ \,{\cal F}^B (\infty) \subseteq {\cal F}^{\widetilde{W}} (\infty) $. 

\smallskip
On the other hand, (\ref{claim})   holds under the condition of case {\it (ii)} in the Proposition, as $\, \langle U \rangle (\cdot)\,$ is then adapted to the filtration    generated by the $n-$dimensional Brownian motion $\, \Gamma (\cdot)\,$; this,  in turn, is independent of  $ \,\widetilde{W}(\cdot)\,$   on the strength of the \textsc{P. L\'evy} Theorem (e.g., \textsc{Karatzas \& Shreve}, Theorem 3.3.16),  since
  $$
  \langle \widetilde{W}, \Gamma_j \rangle (\cdot)\,=\int_0^{\, \cdot} \text{sgn} \big( \widetilde{X} (\theta)  \big)\, \dx \langle B, \Gamma_j \rangle (\theta)\, \equiv \,0\,, \qquad \forall ~~ j=1, \cdots, n\,.   
   $$
The proof of the proposition is complete. \qed

\begin{rem} 
It would be interesting to obtain sufficient conditions for either (\ref{claim})  or (\ref{"pure"}) to hold, which are weaker than those of Proposition \ref{Prop2}. As Example \ref{Counter} shows, however -- and contrary to our own initial guess --  we cannot expect the conclusions of Proposition \ref{Prop2} to remain true for general \textsc{Ocone} martingales. 
\end{rem}

\begin{exmp}
\label{Skew^1} {\it From   Brownian Motion to  Skew Brownian Motion:} Suppose that $\, U(\cdot)\,$ is standard, real valued Brownian motion. Then the conditions of Propositions \ref{Prop1} and \ref{Prop2} are satisfied rather trivially; uniqueness in distribution holds for the skew-\textsc{Tanaka} equation of (\ref{3a}) (equivalently,  of (\ref{3})); and every continuous semimartingale $\,X(\cdot)\,$ that satisfies (\ref{3a}) is of the form
$$
 X (\cdot)  \,=\,  W  (\cdot) + {\, 2 \, \al - 1\, \over \al}\, L^{ X } (\cdot)\qquad \text{with} \qquad W(\cdot)\,:=\int_0^{\, \cdot}  \text{sgn} \big( X(t) \big) \, \dx U(t)\,,
$$
or equivalently 
$$
X (\cdot)  \,=\,  W  (\cdot) + 2\, \big(  2 \, \al - 1\big)\,  \widehat{L}^{ X } (\cdot)
$$
in terms of the symmetric local time as in (\ref{LT}). Of course $\,W(\cdot)\,$ is standard Brownian motion by the \textsc{P. L\'evy} theorem, and the \textsc{Harrison-Shepp} (1981) theory once again characterizes $\,X(\cdot)\,$ as skew Brownian motion with parameter $\, \al\,$.  The processes $\, W(\cdot)\,$ and $\, X(\cdot)\,$ generate the same filtration, which is strictly finer than the filtration generated by the original Brownian motion $\, U(\cdot)=\int_0^{\, \cdot}  \text{sgn} \big( X(t) \big) \, \dx W(t)\,$.
\end{exmp}

\begin{exmp}
\label{Counter} {\it Failure of Uniqueness in Distribution for General \textsc{Ocone}  Martingales:} 
We adapt to our setting a construction from page 131 of \textsc{Dubins et al.}$\,$(1993). We start with a filtered probability space $(\Om, \mathcal{F}, \mathbb{P}), \,\mathbb{F}^B  = \{ \mathcal{F}^B  (t) \}_{0 \le t < \infty}\,$ where $\, B(\cdot)\,$ is standard Brownian motion with $B(0)=0$, and define the adapted, continuous and  strictly increasing process
\begin{equation} \label{eq: ex3}
A({t}) \, :=\, t \cdot {\bf 1}_{\{ t \le 1\}} + \big \{ 1 + \big( u \cdot {\bf 1}_{ \{ B(1) > 0\}} + v \cdot {\bf 1}_{\{B(1) \le 0\}} \big)(t-1) \big\} \cdot {\bf 1}_{ \{ t > 1\}} \, , \quad 0 \le t < \infty
\end{equation}
where $u >0$ and $v>0$ are given real numbers with $u \neq v$, as well as the processes 
\begin{equation} \label{eq: ex1}
X(\cdot) \, :=\, B(A(\cdot)) \, , 
\qquad \Xi(\cdot) \,:=\, -X(\cdot)\,.
\end{equation}
The \textsc{L\'evy} transform
$$
\betab (\cdot) \, :=\,  \int^{\,\cdot}_{0} \text{sgn}  (B(t))\, {\mathrm d} B(t)
$$
of $B(\cdot)$ is a standard Brownian motion adapted to the filtration $   \,\mathbb{F}^{|B|}  = \{ \mathcal{F}^{|B|}  (t) \}_{0 \le t < \infty}\,$, which is strictly coarser than $\,\mathbb{F}^B  \,$; in particular, it can be seen that  $\, \betab (\cdot) \, $ is independent of sgn$(B(1))= 2\, \mathbf{ 1}_{ \{ B(1) >0\} } -1\,$,   and thus of the process $\, A(\cdot)\,$ as well. 

On the other hand, the process $\,X(\cdot)\,$ is a martingale of its natural filtration $  \,\mathbb{F}^X  = \{ \mathcal{F}^B \big( A (t) \big)\}_{0 \le t < \infty}\,$; therefore, so is its ``mirror image" $\, \Xi (\cdot)\,$, and more importantly its \textsc{L\'evy} transform
$$
U (\cdot) \,  :=\,  \int^{\,\cdot}_{0} \text{sgn}  \big(X(t)\big)\, {\mathrm d} X(t) \,=\, \betab \big( A(\cdot)\big) \qquad \text{with} \qquad \langle U  \rangle (\cdot) \, =\,  A(\cdot)\,,
$$
which is thus seen to be an \textsc{Ocone} martingale. Now clearly, both $\, X(\cdot)\,$ and $\, \Xi(\cdot)\,$ satisfy the equation (\ref{3a}) with $\, \alpha = 1/2\,$ driven by $\, U(\cdot)$, so pathwise uniqueness fails for this equation. We also note that the conditions    of Proposition \ref{Prop2} fail  too in this case. 

{\it We claim  that uniqueness in distribution fails as well.} In a  manner similar to the treatment  in \textsc{Dubins et al.} (1993), we shall argue that the distributions of $\,X(\cdot) \,$ and $\, \Xi(\cdot)\,$ at time $\,t \, =\, 2\,$ are different. Now if the random variables 
\[
X(2) \, =\,  B(1+u)\cdot {\bf 1}_{\{B(1) > 0\}} + B(1+v) \cdot {\bf 1}_{\{B(1) \le 0\}}   \qquad \text{and} \qquad \Xi (2) \, =\, - X(2) 
\]
 had the same  probability distributions, that is, if the distribution of the random variable $\, X(2)\,$ were symmetric about the origin, we would have $\, \mathbb{E} [ (X(2))^3] =0\,$. However, let us note the decomposition  
 \[
X(2) \, =\,  B(1) + \big( B(1+u) - B(1) \big) \cdot {\bf 1}_{\{B(1) > 0\}} + \big( B(1+v) - B(1) \big) \cdot {\bf 1}_{\{B(1)\le 0\}}\,,
 \]
 which gives 
 \[ 
 \mathbb{E} \big[ (X(2))^3\big] \,=\, 3\, \mathbb{E} \left[ B(1) \, \big( B(1+u) - B(1) \big)^2 \, {\bf 1}_{\{B(1) > 0\}} \right] + 3\, \mathbb{E} \left[ B(1) \, \big( B(1+v) - B(1) \big)^2 \, {\bf 1}_{\{B(1) \le  0\}} \right]
 \]
 \[
 =\, 3\, \mathbb{E} \left[ \big( B(1) \big)^+\,\right]  \, \big(u-v \big) \, \neq \, 0\,.~~~~~~~~~~~~~~~~~~~~~~~~~~~~~~~~~~~~~~~~~~~~~~~~~~~~~~~~~~~~~~~~~~~~~~~~~ 
 \]
 This contradiction   establishes the claim. 

\end{exmp}

\subsection{The Perturbed Skew-Tanaka Equation is Strongly Solvable}
\label{Str}

The addition of some independent noise can restore pathwise uniqueness, thus also strength, to weak solutions of the  stochastic equation in (\ref{3}) or (\ref{3a}). In the spirit of \textsc{Prokaj} (2013) or \textsc{Fernholz, Ichiba, Karatzas \& Prokaj} (2013), hereafter   referred to as [FIKP], we have the following result.

\begin{prop}
 \label{Prop5}
Suppose that the continuous semimartingale $U(\cdot)$ as in (\ref{1}) satisfies the conditions of Proposition \ref{Prop1}, where now the $\mathbb{F}-$progressively measurable process $\,p(\cdot)\,$ of (\ref{3d}) is locally square-integrable with respect to $\, \langle M \rangle (\cdot)\,;$  and that $$\, V(\cdot) \,=\, N(\cdot) + \Delta (\cdot)\,$$ is another continuous semimartingale, with continuous local martingale part $N(\cdot)$ and finite variation part $\Delta(\cdot)$ which satisfy $N(0)= \Delta(0)=0\,$ and
$$
\langle M, N \rangle (\cdot) \equiv 0\,, \qquad \langle M  \rangle (\cdot) \, =\int_0^{\, \cdot} q (t) \, \dx   \langle N  \rangle ( t)
$$
for some $\mathbb{F}-$progressively measurable process $\,q(\cdot)\,$ with values in a compact interval $\,[0, b]\,$. 

Then pathwise uniqueness holds for the perturbed skew-\textsc{Tanaka} equation 
\begin{equation}
\label{skewProk}
 X(\cdot) \,=\int_0^{\, \cdot}  \text{sgn} \big( X(t) \big)\, \dx U(t)+ V(\cdot) + {\, 2 \, \al - 1\, \over \al}\, L^X (\cdot)\,,
\end{equation}
provided that either \\ (i) $~\,\alpha = 1/2\,$, or that   \\ (ii) $\,U(\cdot)$ and $V(\cdot)$ are independent, standard Brownian motions. In this   case  a weak solution to (\ref{skewProk}) exists, and is thus strong. 
\end{prop}

The claim of case {\it (i)} is proved in Theorem 8.1 of [FIKP], and the claim of case {\it (ii)} in an Appendix, section \ref{PfProp}. 
In case {\it (ii)} of Proposition \ref{Prop5} the equation \eqref{skewProk} can be written equivalently as 
$$
 X(\cdot) \,=  \int_0^{\, \cdot}  \mathbf{ 1}_{ \{ X(t) >0\} } \, \dx W_+(t) \,+  \int_0^{\, \cdot}  \mathbf{ 1}_{ \{ X(t) < 0\} } \, \dx W_-(t) \, + \,{\, 2 \, \al - 1\, \over \al}\, L^X (\cdot) \,.
$$
Here $\, W_\pm (\cdot) :=    V(\cdot)\pm U(\cdot)\,$ are independent Brownian motions with local variance 2; one of them governs the motion of $\, X(\cdot)\,$ during its positive excursions,  the other during the negative ones, whereas these excursions get skewed when $\, \alpha \neq 1/2\,$. 

\section{Proof of Theorem \ref{Thm1}}
\label{Pf}

We shall follow very closely the methodology of \textsc{Prokaj} (2009), with some necessary modifications related to the skewness. The enlargement of the filtered probability space $(\Om, \mathcal{F}, \mathbb{P}), \,\mathbb{F} = \{ \mathcal{F} (t) \}_{0 \le t < \infty}$  is done in terms of a sequence $\, \{ \xi_k\}_{k \in \N}\,$ of independent random variables with common \textsc{Bernoulli} distribution
\begin{equation}
\label{4}
\mathbb{P} \big( \xi_1 = +1 \big) = \al\,, \qquad \mathbb{P} \big( \xi_1 = -1 \big) = 1-\al
\end{equation}
(thus with expectation $\, \mathbb{E} ( \xi_1) =  2 \al -1$), which is independent of $\, \mathcal{F} (\infty) = \sigma \big( \bigcup_{\,0 \le t < \infty} {\cal F}(t) \big)\,$. On the enlarged probability space $\, \big(\widetilde{\Om}, \widetilde{\mathcal{F}}, \widetilde{\mathbb{P}}\big)\,$ we have all the objects of    the original space, so we keep the same notation for them. We denote by 
\begin{equation}
\label{5}
\mathfrak{Z} \,:=\, \big\{ t \ge 0 : S(t) =0 \big\}
\end{equation}
the zero set of the \textsc{Skorokhod} reflection $\, S(\cdot) $ in (\ref{2}), and enumerate as $\, \{ \mathcal{C}_k\}_{k \in \N}\,$ the disjoint components of $\, [0, \infty) \setminus \mathfrak{Z}\,$, that is, the countably-many excursion intervals of the process $\, S(\cdot)\,$ away from the origin. This we do in a measurable manner, so that 
$$
\big\{ t \in \mathcal{C}_k \big\} \in \mathcal{F} (\infty)\,, \qquad \forall ~~ t \ge 0\,, ~ k \in \N\,.
$$
In order to simplify notation, we set 
\begin{equation}
\label{6}
\mathcal{C}_0 \,:=\, \mathfrak{Z} \,, \quad \xi_0 \,:=\,0\,.
\end{equation}

We define now
\begin{equation}
\label{7}
Z(t) \,:=\, \sum_{k \in \N_0} \, \xi_k\, \mathbf{ 1}_{   \mathcal{C}_k  } (t)\,, \qquad \widetilde{\mathcal{F}} (t) := \mathcal{F} (t) \vee \mathcal{F}^Z (t)
\end{equation}
for all $ \, t \in [0, \infty)\,$; this gives the enlarged filtration $ \,\widetilde{\mathbb{F}} = \big\{ \widetilde{\mathcal{F}} (t) \big\}_{0 \le t < \infty}\,$. We posit  the following two claims.

\begin{prop}
\label{Prop3}
The process $\, M(\cdot)\,$ of (\ref{1}) is a continuous local martingale of the enlarged filtration $\, \widetilde{\mathbb{F}}\,$. Consequently, both $\, U(\cdot)\,$ and $\, S(\cdot)\,$ are continuous $ \,\widetilde{\mathbb{F}}-$semimartingales.
\end{prop}

\begin{prop}
\label{Prop4}
In the notation of (\ref{2}) and (\ref{7}), we have 
\begin{equation}
\label{8}
Z(\cdot) \, S(\cdot)\,=\, \int_0^{\, \cdot} Z(t)\, \dx S(t) + \big( 2 \al -1\big)\, L^S (\cdot)\,.
\end{equation}
\end{prop}

\medskip
Taking the claims of these two propositions at face-value for a moment, we can proceed with the proof of Theorem \ref{Thm1} as follows. We define the process 
\begin{equation}
\label{9}
X(\cdot) \,:=\, Z (\cdot) \, S(\cdot)
\end{equation}
and note 
\begin{equation}
\label{10}
Z(\cdot) \,=\, \overline{\text{sgn}} \big( X(\cdot) \big)\,, \qquad \big| X(\cdot) \big| \,=\, S(\cdot)
\end{equation}
thanks to (\ref{6}) and (\ref{7}), as well as 
\begin{equation}
\label{11}
X(\cdot) - \int_0^{\, \cdot} \overline{\text{sgn}} \big( X( t) \big)\,\dx S(t)\,=\, Z (\cdot) \, S(\cdot) - \int_0^{\, \cdot}Z( t)  \,\dx S(t)\,=\, \big( 2 \al -1 \big)\, L^S (\cdot)
\end{equation}
thanks to (\ref{9}), (\ref{8}). In particular, $X(\cdot)\,$ is an $\, \widetilde{\mathbb{F}}-$semimartingale, and we note the property
$$
2\, L^X (\cdot) - L^S (\cdot) \,=\,2\, L^X (\cdot) - L^{|X|} (\cdot) \,=\, \int_0^{\, \cdot} \mathbf{ 1}_{ \{ X(t) =0 \}  }\, \dx X(t)
$$
of its local time at the origin (cf. section 2.1 in \textsc{Ichiba et al.} (2013)). In conjunction with (\ref{11}) and the fact that $X(\cdot)$, $S(\cdot)$, and $Z(\cdot)$ all have the same zero set $\, \mathfrak{Z}\,$ as in (\ref{5}),  (\ref{6}), we get from this last equation 
\begin{equation}
\label{11a}
2\, L^X (\cdot) - L^S (\cdot) \,= \int_0^{\, \cdot} \mathbf{ 1}_{ \{ X(t) =0 \}  }\,\big[ \,  \overline{\text{sgn}} \big( X( t) \big) \,\dx S(t) + \big( 2 \al -1 \big)\, L^S ( t) \, \big]= \big( 2 \al -1 \big)\, L^S (\cdot)\,,~~
\end{equation}
thus 
\begin{equation}
\label{11b}
L^X (\cdot) \,=\, \al\, L^S(\cdot)\,, 
\end{equation}
establishing the second equality in (\ref{3}). Back in (\ref{11}), this leads to 
\begin{equation}
\label{12}
X(\cdot) \,=\, \int_0^{\, \cdot} \overline{\text{sgn}} \big( X( t) \big)\, \big[\, \dx U(t) +  \dx C(t) \big] + \big( 2 \al -1 \big)\, L^S (\cdot)\,,
\end{equation}
where $C(\cdot)$ is the continuous, adapted and increasing process 
$$
C(t) \,:=\, S(t) - U(t) \,=\, \max_{0 \le s \le t} \big( - U(s) \big)\,, \qquad 0 \le t < \infty\,.
$$
From the theory of the \textsc{Skorokhod} reflection problem we know that this process $\,C(\cdot)\,$ is flat off the set $   \{ t \ge 0 : S(t) =0\}  = \mathfrak{Z}\,$, so the skew-\textsc{Tanaka} equation of (\ref{3}) follows now from (\ref{12}), (\ref{11b}). 

The proof of Theorem \ref{Thm1} is complete.   \qed

\medskip
\noindent
{\it Proof of Proposition \ref{Prop3}:} By localization of necessary, it suffices to show that if $M(\cdot)$ is an $\,\mathbb{F}-$martingale, then it is also an $\,\widetilde{\mathbb{F}}-$martingale; that is, for any given $\, 0 < \theta <t<\infty$ and $\, A \in \widetilde{\mathcal{F}} (\theta)\,$ we have 
\begin{equation}
\label{13}
\mathbb{E} \, \big[\, \big( M(t) - M(\theta) \big) \, \mathbf{ 1}_A \, \big] \,=\,0\,.
\end{equation}
It is clear from (\ref{7}) that we need to consider only sets of the form $\, A = B \cap D\,$, where $\, B \in \mathcal{F} (\theta)\,$ and 
\begin{equation}
\label{14}
D \,=\, \bigcap_{j=1}^n \big\{ Z(t_j) = \ep_j \big\} \,=\, \bigcap_{j=1}^n \big\{ \xi_{\,\kappa (t_j)} = \ep_j \big\} 
\end{equation}
for $\, n \in \N\,$, $\, 0 < t_1 < t_2 < \cdots < t_n <\theta <t\,$ and $\, \ep \in \{ -1, 0 , 1\}\,$. Here we have denoted by $\, \kappa (u)\,$ the (random) index of the excursion interval $\, \mathcal{C}_k\,$ to which a given $\, u \in [0, \infty)\,$ belongs. 

For such choices, and because
$$
\mathbb{E} \, \big[\, \big( M(t) - M(\theta) \big) \, \mathbf{ 1}_A \, \big] \,=\,\mathbb{E} \, \big[\, \big( M(t) - M(\theta) \big) \, \mathbf{ 1}_B \cdot \mathbb{E} \, \big( \mathbf{ 1}_D \,|\, \mathcal{F} (\infty) \big) \, \big] \,,
$$
we see that, in order to prove (\ref{13}), it is enough to argue that 
\begin{equation}
\label{15}
\mathbb{E} \, \big( \mathbf{ 1}_D \,|\, \mathcal{F} (\infty) \big) \quad\text{is}~~\mathcal{F} (\theta)-\text{measurable.}
\end{equation}
But the random variables $\, \kappa (t_j)\,$ in (\ref{14}) are measurable with respect to  $\, \mathcal{F} (\infty)\,$, whereas the random variables $\,   \xi_1\,, \, \xi_2\,, \cdots \,$ are independent  of this $\sigma-$algebra. Therefore, we have 
\begin{equation} 
\label{eq: quant}
\mathbb{E} \, \big( \mathbf{ 1}_D \,|\, \mathcal{F} (\infty) \big)= \mathbb{P} \left[\,  \bigcap_{j=1}^n \big\{ \xi_{\,\kappa (t_j)} = \ep_j \big\}  \, \Big| \, \mathcal{F} (\infty)\, \right]=\mathbb{P} \big( \, \xi_{k_1} = \ep_1, \cdots ,  \xi_{k_n} = \ep_n\, \big) \Big|_{k_1 = \kappa (t_1), \cdots, k_n = \kappa (t_n)}.
\end{equation}
For given indices $\, ( k_1, \cdots, k_n)\,$ and $\, ( \ep_1, \cdots, \ep_n)\,$, let us denote by $m$ the number of distinct non-zero indices in $\, ( k_1, \cdots, k_n)\,$,   by $\, \lambda \,$ the number from among those distinct indices of the corresponding $\, \ep_j$'s that are equal to 1, and observe
$$
\mathbb{P} \big( \, \xi_{k_1} = \ep_1, \cdots ,  \xi_{k_n} = \ep_n\, \big) \,=\,0\,, ~\text{if}~ ~ ( \ep_1, \cdots, \ep_n) ~~ \text{contradicts}~~ ( k_1, \cdots, k_n)\,;
$$
\begin{equation}
\label{quant}
~~~~~~~~~~~~~~~\,\,~~~~=\, \al^{\,\lambda} \, \big( 1-\al \big)^{m-\lambda}\,, ~~\text{otherwise}\,.
\end{equation}
\noindent
Here ``$  ( \ep_1, \cdots, \ep_n)\,$ contradicts $\, ( k_1, \cdots, k_n) $" means that we have either \\ $.~ k_i =k_j$ but $\ep_i \neq \ep_j\,$ for some $\, i \neq j\,$; or \\ $.~ k_i=0$ but $\, \ep_i \neq 0\,$, for some $\,i\,$; or \\ $.~ k_i \neq 0$ but $\, \ep_i = 0\,$, for some $\,i\,$.

\medskip
We note now that when $\,{k_1 = \kappa (t_1)\, , \, \cdots \, , \, k_n = \kappa (t_n)}\,$, the value of $\,m\,$ (that is, the number of excursion intervals in $\, [0,s] \setminus \mathfrak{Z} \,$ that contain some $\,t_i$), the value of $\lambda$ (i.e., the number of   such excursion intervals that are positive) and the statement ``$\, ( \ep_1, \cdots, \ep_n)\,$ contradicts $\, ( k_1, \cdots, k_n)\,$", can all be determined on the basis of the trajectory $\, S(u),\, 0 \le u \le \theta\,$; that is, the quantity on the right-hand side of (\ref{eq: quant}) is $\, \mathcal{F}^S(\theta)-$measurable. As a consequence, the property (\ref{15}) holds. \qed

 \medskip
\noindent
{\it Proof of Proposition \ref{Prop4}:} For any $\, \ep \in (0,1)\,$ we define recursively, starting with $\, \tau^\ep_0 :=0\,$, a sequence of stopping times
$$
\tau^\ep_{2 \ell +1 }\,:=\, \inf \big\{ t > \tau^\ep_{2 \ell   }\,:\, S(t) > \ep \big\}\,, \qquad \tau^\ep_{2 \ell +2 }\,:=\, \inf \big\{ t > \tau^\ep_{2 \ell +1  }\,:\, S(t) =0 \big\}
$$
for $\, \ell \in \N_0\,$. We use this sequence to approximate the process $\, Z(\cdot)\,$ of (\ref{7}) by
$$
Z^\ep (t)\,:=\, \sum_{\ell \in \N_0} \, Z(t)\, \mathbf{ 1}_{ \,( \tau^\ep_{2 \ell +1 }, \tau^\ep_{2 \ell +2 }]} (t)\,, \qquad 0 \le t < \infty\,.
$$

Let us note that the resulting process $\, Z^\ep (\cdot)\,$ is constant on each of the indicated intervals; that the sequence of stopping times just defined does not accumulate on any bounded time-interval, on account of the fact that $\, S(\cdot)\,$ has continuous paths; and that the process  $\, Z^\ep (\cdot)\,$ is  of finite first variation over compact intervals. We deduce
\begin{equation}
\label{16}
Z^\ep (T)\, S(T) \,=\, \int_0^{T} Z^\ep ( t)\, \dx S( t) + \int_0^{T} S( t) \, \dx Z^\ep ( t)\,, \qquad 0 \le T <\infty \,. 
\end{equation}
The piecewise-constant process $\, Z^\ep (\cdot)\,$ tends to $\, Z(\cdot)\,$ pointwise as $\, \ep \da 0\,$, and we have 
\begin{equation}
\label{17}
\lim_{\ep \da 0} \, \int_0^{T} Z^\ep ( t)\, \dx S( t)\,=\,  \int_0^{T} Z (t)\, \dx S( t)\,, \quad \text{in probability} 
\end{equation}
for any given $\, T \in [0, \infty)\,$; all the while, $\, |Z^\ep (\cdot)| \le 1\,$. On the other hand, the second integral in (\ref{16}) can be written as 
$$
 \int_0^{T} S( t) \, \dx Z^\ep ( t)\,=  \sum_{ \{ \ell \,: \,\tau^\ep_{2 \ell +1 } < T \}}  S\big( \tau^\ep_{2 \ell +1 } \big) \, Z\big( \tau^\ep_{2 \ell +1 } \big)= \, \ep \sum_{ \{ \ell \,:\, \tau^\ep_{2 \ell +1 } < T \}}    Z\big( \tau^\ep_{2 \ell +1 } \big)
$$
$$
~~~~~~~~~~~~~~= \, \ep \sum_{j =1}^{N(T,\ep)} \xi_{\, \ell_j}\,=\, \ep\, N(T,\ep) \cdot { 1 \over \,N(T,\ep)\,}  \sum_{j =1}^{N(T,\ep)} \xi_{\, \ell_j}\,,
$$
\noindent
where $\, \big\{ \xi_{\ell_j} \big\}_{j=1}^{N(T, \ep)}\,$ is an enumeration of the values $\,Z\big( \tau^\ep_{2 \ell +1 } \big)\,$ and 
$$ N(T, \ep) \,:= \,\# \, \big\{ \ell \,: \,\tau^\ep_{2 \ell +1 } < T \big\}\,$$ is the number of upcrossings of the interval $\, (0, \ep)\,$ that the process $\, S(\cdot)\,$ has completed    by time $T$. From Theorem VI.1.10 in \textsc{Revuz \& Yor} (1999), we have the representation of local time 
$\,
\lim_{\ep \da 0} \, \,\ep \, N (T, \ep) = L^S(T)\,$;  
whereas the strong law of large numbers gives
$$
\lim_{\ep \da 0} \, { 1 \over \,N(T,\ep)\,}  \sum_{j =1}^{N(T,\ep)} \xi_{\, \ell_j}\,=\, \mathbb{E} \big( \xi_1\big)\,.
$$
Back into (\ref{16}) and with the help of (\ref{17}), these considerations give  
$$
Z(T) \, S(T)\,=\, \int_0^{T} Z(t)\, \dx S(t) +\mathbb{E} \big( \xi_1\big)\cdot L^S (T)\,, \qquad 0 \le T <\infty\,,
$$
that is, (\ref{8}). \qed

\section{Conventional Reflection}
\label{CR}

In a similar manner  one can establish the following analogue of Theorem  \ref{Thm1},  which uses the conventional reflection in place of the \textsc{Skorokhod} reflection. 

\begin{thm}
 \label{Thm2}
Fix a constant $\alpha \in (0,1)$. There exists an enlargement $\,\big(\widehat{\Om}, \widehat{\mathcal{F}}, \widehat{\mathbb{P}}\big),\, \widehat{\mathbb{F}} = \{ \widehat{\mathcal{F}} (t) \}_{0 \le t < \infty}$ of the filtered probability space $(\Om, \mathcal{F}, \mathbb{P}), \,\mathbb{F} = \{ \mathcal{F} (t) \}_{0 \le t < \infty}\,$, with a measure-preserving map $\,\pi : \Om \ra \widehat{\Om}\,$, and on this enlarged space a continuous semimartingale $\widehat{X}(\cdot)$    that satisfies
\begin{equation}
\label{3f}
\big|\widehat{X}(\cdot) \big| = \big|U(\cdot) \big|\,,\qquad L^{\widehat{X}} (\cdot) = \al\, L^{|U|}(\cdot)\,, \qquad \widehat{X}(\cdot) \,=\int_0^{\, \cdot} \overline{\text{sgn}} \big( \widehat{X}(t) \big)\, \dx \widehat{U}(t) + {\, 2 \, \al - 1\, \over \al}\, L^{\widehat{X}} (\cdot)\,. 
\end{equation}
Here 
\begin{equation}
\label{18}
\widehat{U}(\cdot)\,:=\, \int_0^{\, \cdot} \overline{\text{sgn}} \big( U(t) \big)\, \dx  U (t)
\end{equation}
is the \textsc{L\'evy}   transform  of the semimartingale $\, U(\cdot)\,$, and the classical reflection $\,   R(\cdot)= | U (\cdot)|\,$ of $\, U (\cdot)\,$ coincides with the \textsc{Skorokhod} reflection  of the process $\, \widehat{U}(\cdot)\,$ in (\ref{18}), namely  
$$\,\widehat{S}(t)\, :=\, \widehat{U}(t) + \max_{0 \le s \le t} \big( - \widehat{U}(s) \big)\,, ~~~~~~ 0 \le t < \infty\,.$$
\end{thm}

Indeed, most of the argument of the   proof in section \ref{Pf} goes through verbatim, with $\, S(\cdot), \, X(\cdot)\,$ replaced here by $\, R(\cdot), \, \widehat{X}(\cdot)\,$, up to and including the display (\ref{11b}). But now we have
\begin{equation}
\label{18a}
R(\cdot)= | U(\cdot)| = \int_0^{\, \cdot}\overline{\text{sgn}} \big( U(t) \big)\, \dx  U (t) + L^{|U|} (\cdot)\, =\, \widehat{U}(\cdot) + L^{R} (\cdot) 
\end{equation}
from the \textsc{It\^ o-Tanaka} formula, so (\ref{12}) is replaced by
$$
\widehat{X}(\cdot) \,=\int_0^{\, \cdot}\overline{\text{sgn}} \big( \widehat{X}(t) \big)\,\big[ \dx  \widehat{U} (t) + \dx L^{R} (  t) \big] + \big( 2 \al -1 \big)  L^R (\cdot)\,. 
$$
The property $\,L^{\widehat{X}} (\cdot) = \al\, L^R(\cdot)\,$ is established   exactly as in (\ref{11b}), so the stochastic integral equation in (\ref{3f}) follows from this last display. On the other hand, since the local time $\, L^R(\cdot)\,$ grows only on the set $\, \{ t \ge 0: R(t) =0\} = \{ t \ge 0:  \widehat{X} (t)=0\}$, the equality of the first and last terms in (\ref{18a}) identifies $\, R(\cdot)\,$ as the \textsc{Skorokhod} reflection $\,\widehat{S}(\cdot)\,$ of the \textsc{L\'evy} transform $\,\widehat{U}(\cdot)$, as claimed in the last sentence of  Theorem \ref{Thm2}. It is well-known (see, for instance,  \textsc{Chaleyat-Maurel \& Yor } (1978)) that the processes $\, | U(\cdot)|\,$ and $\,\widehat{U}(\cdot)\,$ generate the same filtration. 

\begin{rem}
Let us note that the stochastic integral equation in (\ref{3f}) can   always be written in the more conventional form 
\begin{equation}
\label{3g}
\widehat{X}(\cdot) \,=\int_0^{\, \cdot}  \text{sgn}  \big( \widehat{X}(t) \big)\, \dx \widehat{U}(t) + {\, 2 \, \al - 1\, \over \al}\, L^{\widehat{X}} (\cdot)\,,
\end{equation}
{\it without any additional conditions on} $U(\cdot)$. This is because the analogue $\,\int_0^{\, \cdot}  \mathbf{ 1}_{ \{ \widehat{X}(t) =0\} }\, \dx \widehat{U}(t)\equiv 0\,$ of the property in (\ref{3b}) is now satisfied trivially, on account of 
(\ref{18}). 
\end{rem}

\begin{exmp}
\label{Skew^2} {\it From One Skew Brownian Motion to Another:} Suppose that $\, U(\cdot)\,$ is a skew Brownian motion with parameter $\, \ga  \in (0,1)  $, i.e., 
$$
U(\cdot ) \,=\, B(\cdot) + { \, 2 \,\ga -1\, \over \ga}\, L^U (\cdot)
$$
for some standard, real-valued Brownian motion $B(\cdot)$. We have in this case $\,\int_0^{\infty} \mathbf{ 1}_{ \{ U(t) =0 \}  }\, \dx t=0\,$ as well as the local time property 
$$\, 
2\, L^U (\cdot) - L^{|U|} (\cdot) = \int_0^{\, \cdot} \mathbf{ 1}_{ \{ U(t) =0 \}  }\, \dx U(t)= \frac{\,2 \,\ga - 1\,}{ \ga}\,  L^U(\cdot)\,,
$$ 
thus $\, L^U (\cdot) = \ga \,L^{|U|} (\cdot)\,$ and therefore
$ 
R(\cdot)\,=\, \big| U(\cdot) \big| \,=\, \int_0^{\, \cdot} \overline{\text{sgn}} \big( U(t) \big) \, \dx U(t) + L^{|U|} (\cdot)\,=\, W(\cdot) + L^{|U|} (\cdot)\,$.
Here we have denoted the \textsc{L\'evy} transform of (\ref{18}) as  
$$
W(\cdot)\,:=\, \widehat{U}(\cdot)\,=\int_0^{\, \cdot} \overline{\text{sgn}} \big( U(t) \big)  \left( \dx  B (t)+ { \, 2 \,\ga -1\, \over \ga}\, \dx L^U ( t) \right) \,=\int_0^{\, \cdot}  \text{sgn}  \big( U(t) \big) \, \dx B(t)\,,
$$
and observed that it is another standard Brownian motion. Thus,  the stochastic integral equation of (\ref{3g}) becomes
$$
\widehat{X}(\cdot) \,=\int_0^{\, \cdot}  \text{sgn}  \big( \widehat{X}(t) \big)\, \dx W (t) + {\, 2 \, \al - 1\, \over \al}\, L^{\widehat{X}} (\cdot)\,=\, \widehat{W} (\cdot) + {\, 2 \, \al - 1\, \over \al}\, L^{\widehat{X}} (\cdot)
$$
with $\, \widehat{W} (\cdot)=\int_0^{\, \cdot}  \text{sgn}  \big( \widehat{X}(t) \big)\, \dx W (t)\,$ yet another standard Brownian motion. 

\smallskip
The \textsc{Harrison-Shepp} (1981) theory characterizes now   $\widehat{X} (\cdot)$  
as   skew Brownian motion  with skewness  parameter $\al\,$. The processes $\widehat{X} (\cdot)$ and $\widehat{W} (\cdot)$ generate the same filtration, as do the processes $$\, \widehat{U}(\cdot) \,=\int_0^{\, \cdot}  \text{sgn}  \big( \widehat{X}(t) \big)\, \dx \widehat{W} (t) =W (\cdot)\, \qquad \text{and } \qquad \,R(\cdot)= | U(\cdot)| \, ;$$ and the first filtration is finer than the second. 
\end{exmp}

\subsection{Skew Bessel Processes}
\label{BESd} 

In this subsection suppose that $\,U^{2}(\cdot)\,$ is a squared \textsc{Bessel} process with dimension $\,\delta \in (1, 2)\,$, i.e., $\,U^{2}(\cdot)\,$ is the unique strong solution of the equation 
$$
U^{2}(t) \, =\,  \delta\,  t + 2 \int^{t}_{0} \sqrt{U^{2}(t)} \,{\mathrm d} B(t) \, , \quad 0 \le t < \infty \, 
$$
for some standard, real-valued Brownian motion $\,B(\cdot)\,$. When $\, \delta \in (1, 2)\,$, the   square root $\,  R(\cdot) := \lvert U(\cdot)\rvert \ge 0\,$ of this process is a semimartingale that keeps visiting  the origin almost surely, and can be decomposed as 
\begin{equation} \label{eq: BESd}
R(\cdot) \, =\,  \int^{\cdot}_{0}\frac{\, \delta - 1\, }{2 \, R(t) \, } \cdot {\bf 1}_{\{R(t) \neq 0\}} {\mathrm d} t + B(\cdot) \, \quad \text{ with } \quad  L^{R}(\cdot) \equiv 0\, , \quad \int^{\cdot}_{0} {\bf 1}_{\{R(t) \, =\, 0\}} {\mathrm d} t \, \equiv\, 0 \, . 
\end{equation}
For the study of the stochastic differential equation (\ref{eq: BESd}) with $\, \delta \in (1, 2) \,$ see, for example, \textsc{Cherny} (2000). 

Given $\,\alpha \in (0, 1) \,$, following again the argument of the proof in section 3 through verbatim, with $\,S(\cdot)\,$, $\,X(\cdot)\,$ replaced respectively by $\,R(\cdot)\,$, $\, \widehat{X}(\cdot)\,$, we unfold the nonnegative \textsc{Bessel} process $\,R(\cdot) \,$ to obtain  
\begin{equation} \label{eq: sBESd}
\widehat{X}(\cdot) \, =\, Z(\cdot) R(\cdot) \, =\,  \int^{\cdot}_{0} Z(t) {\mathrm d} R(t) + (2 \alpha - 1) L^{R}(\cdot) 
\, =\, \int^{\cdot}_{0}\frac{\, \delta - 1\, }{2 \, \widehat{X} (t) \, } \cdot {\bf 1}_{\{ \widehat{X}(t) \neq 0\}} {\mathrm d} t + \widehat{{\bm \beta}} (\cdot) \,,
\end{equation}
with $\, Z(\cdot) \:= \text{sgn}  ( \widehat{X}(\cdot)) \,$ and with $\, \widehat{{\bm \beta}}(\cdot) :=\int^{\cdot}_{0} Z(t) {\mathrm d} B(t) \,$ another standard Brownian motion on an extended probability space, as a consequence of Theorem 4.1 and of the properties in (\ref{eq: BESd}). We note that the semimartingale  $\, \widehat{X}(\cdot)\,$   does not accumulate local time  at the origin, because of $\, L^{R}(\cdot) \equiv 0\,$. 

\smallskip
We claim that the process $\, \widehat{X}(\cdot) \,$ constructed here in (\ref{eq: sBESd}) is the {\it $\,\delta\,$-dimensional  skew Bessel process with skewness parameter $\,\alpha \,$}. This process was  introduced and studied in \textsc{Blei} (2012).  

\smallskip
Indeed, let us consider the functions $\,g(x) \, :=\, \lvert x \rvert^{2-\delta} / (2-\delta)\,$ and $\,G(x) \, :=\, \text{sgn} (x) \cdot g(x)  \,$ for $\, x \in \mathbb R \,$, and examine $\, g( \widehat{X}(\cdot))\,$ and $\, G( \widehat{X}(\cdot))\,$. This scaling is a right choice to measure the boundary behavior of $\, \widehat{X}(\cdot) \,$ around the origin. By substituting $\,q = 2-\delta\,$, $\,p = (2-\delta)\, / \, (1-\delta)\,$, $\,\nu = - 1\, / \, 2 \,$ in Proposition XI.1.11 of \textsc{Revuz \& Yor} (2005), we find there exists a (nonnegative) one-dimensional \textsc{Bessel} process $\, {\bm \rho} (\cdot) \,$ on the same probability space such that $\, {\bm \rho} (0) \, =\,  (2-\delta)^{\delta-1} g( \widehat{X}(0))\,$ and 
\[
g( \widehat{X}(t)) \, =\, \frac{1}{\, 2-\delta\, } \, \big|  \widehat{X}(t)\big|^{2-\delta} \, =\, \frac{1}{\, (2-\delta)^{\delta - 1}} \, {\bm \rho} \big(\Lambda (t)\big) \, , \qquad 0 \le t < \infty\,,
\]
where 
\[
\Lambda (t)  \, :=\, \inf\{ s \ge 0: K(s)  \ge t \}\, , \quad  K(s) \, :=\, \int^{s}_{0} \big({\bm \rho}(u)\big)^{ \frac{2\delta - 2}{2 - \delta}} {\mathrm d} u  \,  , 
\] 
that is, $\,g( \widehat{X}(\cdot))\,$ is a time-changed, conventionally reflected Brownian motion with the stochastic clock $\, \Lambda (\cdot)\,$. Thus the local time of $\,g (\widehat{X}(\cdot))\,$ accumulates at the origin with this clock $\, \Lambda (\cdot)\,$. 

In the same manner as in the  construction of $\, Z(\cdot) R(\cdot) \,$ in Theorem 4.1, we obtain here 
\[
G( \widehat{X}(T)) \, =\,  \text{sgn} ( \widehat{X}(T)) g( \widehat{X}(T)) \, =\,  \int^{T}_{0} \text{sgn}  ( \widehat{X}(t) ) {\mathrm d} \big( g( \widehat{X}(t))\big) +  (2 \alpha - 1) \,  L^{g(\widehat{X})}(T) \, 
\]
as well as 
\begin{equation} \label{eq: skBESlt}
L^{G( \widehat{X})}(\cdot) - L^{-G( \widehat{X})} (\cdot) \, =\,  (2 \alpha - 1) \big( L^{G (\widehat{X})}(\cdot) + L^{-G( \widehat{X})}(\cdot)\big)  
\end{equation} 
and
\[
(1-\alpha) \, L^{G( \widehat{X})}(\cdot)  \, =\,  \alpha \, L^{-G( \widehat{X})}(\cdot) \, , \quad L^{g( \widehat{X})}(\cdot) 
\, =\,  \frac{\, 1\, }{\, 2\, } \big( L^{G (\widehat{X})}(\cdot) + L^{-G( \widehat{X})}(\cdot)\big) \, . 
\]
in the notation of (\ref{LT}).  
From these relationships (\ref{eq: skBESlt}), and on the strength of Theorem 2.22 of \textsc{Blei} (2012), we identify the process of (\ref{eq: sBESd}) as the $\,\delta\,$-dimensional skew \textsc{Bessel} process. Here the process $\,G ( \widehat{X}(\cdot)) \,$ and its local time $\,L^{G( \widehat{X})}(\cdot)\,$ correspond to $\, Y(\cdot)\,$ and $\,L^{X}_{m}(\cdot)\,$, respectively,  in the notation of \textsc{Blei} (2012). 

\smallskip
For various properties and representations of this process, we refer the study of \textsc{Blei} (2012), in particular, Remark 2.26 there.

\section{An Application: Two Diffusive Particles with   Asymmetric Collisions}
\label{Appl}

In the paper [FIKP], the authors construct a planar continuous semimartingale  $\,\mathcal  X (\cdot)= (X_{1}(\cdot), X_{2}(\cdot))\,$ with dynamics 
\begin{equation} 
\label{eq: 2D1}
{\mathrm d} X_{1}(t) = \big( g {\bf 1}_{\{X_{1}(t) \le X_{2}(t)\}} - h {\bf 1}_{\{X_{1}(t) > X_{2}(t)\}} \big) {\mathrm d} t + \big( \rho {\bf 1}_{\{X_{1}(t) > X_{2}(t)\}} + \sigma {\bf 1}_{\{X_{1}(t) \le X_{2}(t)\}}\big) {\mathrm d} B_{1}(t)\,   , ~~~
\end{equation}
\begin{equation} 
\label{eq: 2D2}
{\mathrm d} X_{2}(t) = \big( g {\bf 1}_{\{X_{1}(t) > X_{2}(t)\}} - h {\bf 1}_{\{X_{1}(t) \le X_{2}(t)\}} \big) {\mathrm d} t + \big( \rho {\bf 1}_{\{X_{1}(t) \le X_{2}(t)\}} + \sigma {\bf 1}_{\{X_{1}(t) > X_{2}(t)\}}\big) {\mathrm d} B_{2}(t)  \, , ~~~
\end{equation} 

\medskip
\noindent 
for arbitrary real constants $\,g, \, h$ and $ \, \rho > 0\,, \, \sigma > 0\,$ with $ \, \rho^2 + \sigma^2=1 $. They show  that, for an arbitrary initial condition $\,  (X_{1}(0), X_{2}(0)) =(x_1, x_2) \in \R^2\,$  and with $\, (B_1 (\cdot), \, B_2(\cdot))\,$ a planar Brownian motion, the system of (\ref{eq: 2D1}), (\ref{eq: 2D2}) has a pathwise unique, strong solution. 

This is a model for two ``competing" Brownian particles, with diffusive motions whose drift and dispersion characteristics are assigned according to their ranks.

  \medskip
\noindent 
$\bullet~$ In another recent paper  \textsc{Fernholz, Ichiba \& Karatzas} (2013), hereafter   referred to as [FIK], a planar continuous semimartingale   $\, \widetilde{\mathcal  X } (\cdot) \, =\, ( \widetilde{X}_{1}(\cdot), \widetilde{X}_{2}(\cdot)) \,$ is constructed according to the dynamics 
\[
{\mathrm d} \widetilde{X}_{1}(t) \, =\, \big( g {\bf 1}_{\{ \widetilde{X}_{1}(t) \le \widetilde{X}_{2}(t)\}} - h {\bf 1}_{\{ \widetilde{X}_{1}(t) > \widetilde{X}_{2}(t)\}} \big) {\mathrm d} t + \big( \rho {\bf 1}_{\{ \widetilde{X}_{1}(t) > \widetilde{X}_{2}(t)\}} + \sigma {\bf 1}_{\{ \widetilde{X}_{1}(t) \le \widetilde{X}_{2}(t)\}}\big) {\mathrm d} \widetilde{B}_{1}(t)  \, 
\]
\begin{equation} \label{eq: 2Dskew1}
+ \,\frac{1-\zeta_{1}}{2}\, {\mathrm d} L^{ \widetilde{X}_{1} - \widetilde{X}_{2}}(t) + \frac{1-\eta_{1}}{2}\, {\mathrm d} L^{ \widetilde{X}_{2}- \widetilde{X}_{1}}(t) \, , 
\end{equation}
\[
{\mathrm d} \widetilde{X}_{2}(t) \, =\, \big( g {\bf 1}_{\{ \widetilde{X}_{1}(t) > \widetilde{X}_{2}(t)\}} - h {\bf 1}_{\{ \widetilde{X}_{1}(t) \le \widetilde{X}_{2}(t)\}} \big) {\mathrm d} t + \big( \rho {\bf 1}_{\{ \widetilde{X}_{1}(t) \le \widetilde{X}_{2}(t)\}} + \sigma {\bf 1}_{\{ \widetilde{X}_{1}(t) > \widetilde{X}_{2}(t)\}}\big) {\mathrm d} \widetilde{B}_{2}(t)  \,
\]
\begin{equation} \label{eq: 2Dskew2}
+\, \frac{1-\zeta_{2}}{2}\, {\mathrm d}  L^{ \widetilde{X}_{1} - \widetilde{X}_{2}}(t) + \frac{1-\eta_{2}}{2} \,{\mathrm d} L^{ \widetilde{X}_{2}- \widetilde{X}_{1}}(t) \, , 
\end{equation}

\smallskip
\noindent
Here again $\,g, \, h$  are arbitrary real constants, $ \, \rho > 0\,$ and $ \, \sigma > 0\,$ satisfy $ \, \rho^2 + \sigma^2=1 \,$,   whereas $\,\zeta_{i}, \eta_{i}\,$ are real constants   satisfying 
\[
0 \le \alpha \, :=\, \frac{\eta}{\, \eta + \zeta\, } \le 1\, , \quad \zeta \, :=\,  1+ \frac{\, \zeta_{1} - \zeta_{2}\, }{2} \, , \quad \eta \, :=\, 1 - \frac{\, \eta_{1} - \eta_{2}\, }{2} \, , \quad \zeta + \eta \neq 0 \, . 
\]
This new system is a version of the previous competing Brownian particle system, but now with {\it elastic and asymmetric collisions} whose effect is modeled by the local time terms $\, L^{ \widetilde{X}_{2}- \widetilde{X}_{1}}(\cdot)\,$  and $\, L^{ \widetilde{X}_{2}- \widetilde{X}_{1}}(\cdot)\,$. Every time the two particles collide, their trajectories feel a ``drag" proportional to these local time terms, whose presence makes the analysis of the system (\ref{eq: 2Dskew1}), (\ref{eq: 2Dskew2}) considerable more involved than that of (\ref{eq: 2D1}), (\ref{eq: 2D2}).

It is shown  in [FIK] under the above conditions that, for an arbitrary initial condition $\,   (\widetilde{X}_{1}(0), \widetilde{X}_{2}(0) ) =(x_1, x_2) \in \R^2\,$, and with $\, (\widetilde{B}_1 (\cdot), \, \widetilde{B}_2(\cdot))\,$ a planar Brownian motion, the system of (\ref{eq: 2Dskew1}), (\ref{eq: 2Dskew2}) has a pathwise unique, strong solution.

\medskip
\noindent 
$\bullet~$ We shall show how to use the unfolding of Theorem \ref{Thm1}, in order to construct the planar process $\, \widetilde{\mathcal  X } (\cdot) = ( \widetilde{X}_{1}(\cdot), \widetilde{X}_{2}(\cdot)) \,$ of (\ref{eq: 2Dskew1}), (\ref{eq: 2Dskew2}) with skew-elastic collisions, starting from the planar diffusion   $\,\mathcal  X (\cdot)= (X_{1}(\cdot),$ $ X_{2}(\cdot))\,$  of (\ref{eq: 2D1}), (\ref{eq: 2D2}). For simplicity, we shall take the initial condition $\, (x_1, x_2) =(0,0)\,$ from now on.

\begin{thm}
 \label{Prop7}
Suppose we are given a planar continuous semimartingale  $\,\mathcal  X (\cdot)=(X_{1}(\cdot), X_{2}(\cdot))\,$ that satisfies the system   of (\ref{eq: 2D1}), (\ref{eq: 2D2})  on some filtered  probability space $\,\big(\Omega, \mathcal F, \mathbb P  \big),\,  \mathbb{F}  =   \{  \mathcal{F}  (t) \}_{0 \le t < \infty}\,$ with a    planar Brownian motion $\, (B_1 (\cdot), \, B_2(\cdot))\,$.

\smallskip
There exists then an enlargement $\, \big(\widetilde{\Omega}, \widetilde{\mathcal F}, \widetilde{\mathbb P} \big),\, \widetilde{\mathbb{F}} = \{ \widetilde{\mathcal{F}} (t) \}_{0 \le t < \infty} \,$ of this filtered  probability space,  with a planar Brownian motion $\, \big(\widetilde{B}_1 (\cdot), \, \widetilde{B}_2(\cdot)\big)\,$, and on it a planar continuous semimartingale $\,\widetilde{\mathcal  X } (\cdot) = (\widetilde{X}_{1}(\cdot), \widetilde{X}_{2}(\cdot))\,$ that satisfies the system of (\ref{eq: 2Dskew1}), (\ref{eq: 2Dskew2}) with skew-elastic collisions,  as well as  
\[
\big(X_{1}( t) - X_{2}( t)\big) + \sup_{0 \le s \le  t } \big(X_{1}(s) - X_{2}(s)\big)^{+} 
\, =\, \big| \widetilde{X}_{1}( t) - \widetilde{X}_{2}( t) \big| \,, \qquad 0 \le t < \infty\, . 
\]
\end{thm}

\smallskip
In other words, the size of the gap between the new processes $\, \widetilde{X}_{1}(\cdot)\,,\, \widetilde{X}_{2}(\cdot)\,$ coincides with the \textsc{Skorokhod} reflection of the difference  $\, X_{1}(\cdot)- X_{2}(\cdot)\,$ of the original processes about the origin. We devote the remainder of this section to the proof of this result.


\subsection{Reduction to symmetric local times}

First, some preparatory steps. We define  the averages $\, \overline{\zeta} \, :=\, (\zeta_{1}+ \zeta_{2}) \, / \, 2\,$, $\,~ \overline{\eta} \, :=\, (\eta_{1} + \eta_{2}) \, / \, 2\,$, and introduce yet another parameter   
\begin{equation}
\label{beta}
\beta \, :=\, \alpha \cdot \frac{\, \zeta_{1} + \zeta_{2}\, }{2} + (1-\alpha) \cdot \frac{\, \eta_{1} + \eta_{2}\, }{2} \, =\, \alpha \, \overline{\zeta} + (1-\alpha) \, \overline{\eta} \, . 
\end{equation}
For notational simplicity we shall write all the processes related to the skew collisions with a tilde,   e.g., $\, \widetilde{Y}(\cdot) \, :=\, \widetilde{X}_{1}(\cdot) - \widetilde{X}_{2}(\cdot) \,$. 
From the relation between the {\it right} local time $\,L^{ \widetilde{Y}} (\cdot)\,$  and the {\it symmetric local time}  $ \, \widehat{L}^{\, \widetilde{Y}}(\cdot)\,$ as in (\ref{LT}),  we obtain the relations

\begin{equation} 
\label{eq: 2Dlt}
\zeta L^{ \widetilde{Y}}(\cdot) \, =\, \eta L^{- \widetilde{Y}}(\cdot)\, , \quad L^{ \widetilde{Y}}(\cdot) \, =\, 2 \, \alpha \,\widehat{L}^{ \, \widetilde{Y}}(\cdot) \, , \quad L^{ \widetilde{Y}}_{-}(\cdot) \, :=\, L^{ - \widetilde{Y}}(\cdot) \, =\, 2\, (1-\alpha) \,\widehat{L}^{\, \widetilde{Y}}(\cdot) \,  
\end{equation}

\medskip
\noindent
as   in  [FIK]. This way, the system (\ref{eq: 2Dskew1})-(\ref{eq: 2Dskew2}) can be re-cast as 
\[
{\mathrm d} \widetilde{X}_{1}(t) \, =\, \big( g {\bf 1}_{\{ \widetilde{X}_{1}(t) \le \widetilde{X}_{2}(t)\}} - h {\bf 1}_{\{ \widetilde{X}_{1}(t) > \widetilde{X}_{2}(t)\}} \big) {\mathrm d} t + \big( \rho {\bf 1}_{\{ \widetilde{X}_{1}(t) > \widetilde{X}_{2}(t)\}} + \sigma {\bf 1}_{\{ \widetilde{X}_{1}(t) \le \widetilde{X}_{2}(t)\}}\big) \,{\mathrm d} \widetilde{B}_{1}(t)  \, 
\]
\begin{equation} \label{eq: 2Dskew3}
{} +  (2\alpha - \beta) \,{\mathrm d} \widehat{L}^{ \,\widetilde{Y}}(t) \, , 
\end{equation}
\[
{\mathrm d} \widetilde{X}_{2}(t) \, =\, \big( g {\bf 1}_{\{ \widetilde{X}_{1}(t) > \widetilde{X}_{2}(t)\}} - h {\bf 1}_{\{ \widetilde{X}_{1}(t) \le \widetilde{X}_{2}(t)\}} \big) {\mathrm d} t + \big( \rho {\bf 1}_{\{ \widetilde{X}_{1}(t) \le \widetilde{X}_{2}(t)\}} + \sigma {\bf 1}_{\{ \widetilde{X}_{1}(t) > \widetilde{X}_{2}(t)\}}\big) {\mathrm d} \widetilde{B}_{2}(t)  \,
\]
\begin{equation} \label{eq: 2Dskew4}
{} +  ( 2- 2\alpha - \beta) {\mathrm d} \widehat{L}^{ \, \widetilde{Y}}(t) \, . 
\end{equation}
We shall construct the system (\ref{eq: 2Dskew3})-(\ref{eq: 2Dskew4}) first, and then   obtain from it the system (\ref{eq: 2Dskew1})-(\ref{eq: 2Dskew2}).  

\subsection{Proof of Theorem \ref{Prop7}}

By applying a \textsc{Girsanov}  change of measure twice, we can remove the drifts from both of the systems (\ref{eq: 2D1})-(\ref{eq: 2D2}) and (\ref{eq: 2Dskew3})-(\ref{eq: 2Dskew4}). Then, in the following, let us construct the two-dimensional Brownian motion with rank-based dispersions and skew-elastic collisions 
\begin{equation} 
\label{eq: noDskew}
\begin{split}
{\mathrm d} \widetilde{X}_{1}(t) \, &=\,  \big( \rho {\bf 1}_{\{ \widetilde{X}_{1}(t) > \widetilde{X}_{2}(t)\}} + \sigma {\bf 1}_{\{ \widetilde{X}_{1}(t) \le \widetilde{X}_{2}(t)\}}\big) {\mathrm d} \widetilde{B}_{1}(t)  \, 
+ (2\alpha - \beta)\, {\mathrm d} \widehat{L}^{\, \widetilde{Y}}(t) \, , \\
{\mathrm d} \widetilde{X}_{2}(t) \, &=\, \big( \rho {\bf 1}_{\{ \widetilde{X}_{1}(t) \le \widetilde{X}_{2}(t)\}} + \sigma {\bf 1}_{\{ \widetilde{X}_{1}(t) > \widetilde{X}_{2}(t)\}}\big) {\mathrm d} \widetilde{B}_{2}(t) 
+  ( 2- 2\alpha - \beta) \,{\mathrm d} \widehat{L}^{\, \widetilde{Y}}(t) \,  
\end{split}
\end{equation}
from the   solution $\,((X_{1}(\cdot), X_{2}(\cdot)), (B_{1}(\cdot), B_{2}(\cdot)))\,$ 
of the system 
\begin{equation} 
\label{eq: noD}
\begin{split}
{\mathrm d} X_{1}(t) \, &=\, \big( \rho {\bf 1}_{\{X_{1}(t) > X_{2}(t)\}} + \sigma {\bf 1}_{\{X_{1}(t) \le X_{2}(t)\}}\big) \,{\mathrm d} B_{1}(t)  \, , \\
{\mathrm d} X_{2}(t) \, &=\, \big( \rho {\bf 1}_{\{X_{1}(t) \le X_{2}(t)\}} + \sigma {\bf 1}_{\{X_{1}(t) > X_{2}(t)\}}\big) \,{\mathrm d} B_{2}(t)  \, , 
\end{split}
\end{equation}
 which is known from [FIKP] to be strongly solvable. 
 Since there is no drift in these last equations, the difference $\,Y(\cdot) \, :=\, X_{1}(\cdot) - X_{2}(\cdot)\,$ between the two components of the system (\ref{eq: noD})  is given by the real-valued  Brownian motion
\begin{equation} \label{eq: Y}
  Y(\cdot) \, =\,    W(\cdot) \, :=\, \rho W_{1}(\cdot) + \sigma W_{2}(\cdot) \, .  
\end{equation}
Here
\[
 W_{1}(\cdot) \, :=\, \int_0^{\,\cdot} {\bf 1}_{\{X_{1}(t) > X_{2}(t) \}} {\mathrm d} B_{1}(t) - \int_0^{\,\cdot} {\bf 1}_{\{X_{1}(t) \le X_{2}(t)\}} {\mathrm d} B_{2}(t) \, , 
\]
\[
 W_{2}(t) \, :=\, \int_0^{\,\cdot} {\bf 1}_{\{X_{1}(t) \le X_{2}(t) \}} {\mathrm d} B_{1}(t) - \int_0^{\,\cdot} {\bf 1}_{\{X_{1}(t) > X_{2}(t)\}} {\mathrm d} B_{2}(t)  
\]

\medskip
\noindent
are independent Brownian motions. As in  [FIKP],  let us recall also the Brownian motion
 \[
V(\cdot)  \, :=\, \rho V_{1}(\cdot) + \sigma V_{2}(\cdot) \, , 
 \]
 where again
\[
 V_{1}(\cdot) \, :=\,  \int_0^{\,\cdot}{\bf 1}_{\{X_{1}(t) > X_{2}(t) \}} {\mathrm d} B_{1}(t) +  \int_0^{\,\cdot}{\bf 1}_{\{X_{1}(t) \le X_{2}(t)\}} {\mathrm d} B_{2}(t) \, , 
\]
\[
 V_{2}(\cdot) \, :=\,  \int_0^{\,\cdot}{\bf 1}_{\{X_{1}(t) \le X_{2}(t) \}} {\mathrm d} B_{1}(t) +  \int_0^{\,\cdot}{\bf 1}_{\{X_{1}(t) > X_{2}(t)\}} {\mathrm d} B_{2}(t)  
\]

 \smallskip
\noindent
 are independent Brownian motions. For a given number $ \alpha \in (0, 1) $,  there exists by Theorem \ref{Thm1}  an adapted, continuous process $\, \widetilde{Y}(\cdot) \,$ which satisfies  
$$ 
Y(t) + \sup_{0 \le s \le t} (- Y(s))^{+} \, =\, 
\big| \widetilde{Y}(t)\big| \, , \qquad  0 \le t < \infty\,
$$
as well as 
\begin{equation} 
\label{eq: 2Dunfold}
  \widetilde{Y}(\cdot) \, = \int_0^{\, \cdot}  \overline{\text{sgn}} \big(\widetilde{Y}(t)\big)\, {\mathrm d} Y(t) +  \frac{\,2\alpha - 1\,}{  \alpha} \,   L^{ \widetilde{Y}}(\cdot)  \, = \int_0^{\, \cdot} {\text{sgn}} \big(\widetilde{Y}(t)\big) \,{\mathrm d} W(t) + 2(2\alpha - 1) \,  \widehat{L}^{ \,\widetilde{Y}}(\cdot)\, ,
\end{equation}
where the last equality follows from   
Proposition \ref{Prop1} and (\ref{eq: 2Dlt}). Thus, the ``unfolded process" $\,\widetilde{Y}(\cdot)\,$ is  a skew Brownian motion, with   skewness parameter $\,\alpha\,$.  

Now  let us define the new planar Brownian motion $\, \big(\widetilde{B}_{1}(\cdot), \widetilde{B}_{2}(\cdot)\big)\,$ as 
\[
{\mathrm d}\widetilde{B}_{1}(\cdot) \, :=\,  \big( {\bf 1}_{\{ Y(\cdot) > 0 , \widetilde{Y}(\cdot) > 0\}} - {\bf 1}_{\{ Y(\cdot) \le 0 , \widetilde{Y}(\cdot) \le 0\}} \big){\mathrm d}B_{1}(\cdot) 
+ \big( {\bf 1}_{\{ Y(\cdot) > 0 , \widetilde{Y}(\cdot) \le 0\}} - {\bf 1}_{\{ Y(\cdot) \le 0 , \widetilde{Y}(\cdot) > 0\}} \big){\mathrm d}B_{2}(\cdot)\, , 
\]
\[
{\mathrm d}\widetilde{B}_{2}(\cdot) \, :=\,  \big( {\bf 1}_{\{ Y(\cdot) > 0 , \widetilde{Y}(\cdot) \le 0\}} - {\bf 1}_{\{ Y(\cdot) \le 0 , \widetilde{Y}(\cdot) > 0\}} \big){\mathrm d}B_{1}(\cdot) 
+ \big( {\bf 1}_{\{ Y(\cdot) > 0 , \widetilde{Y}(\cdot) > 0\}} - {\bf 1}_{\{ Y(\cdot) \le 0 , \widetilde{Y}(\cdot) \le 0\}} \big){\mathrm d}B_{2}(\cdot) \, , 
\]

\medskip
\noindent
and, with the number   $\, \beta \in \mathbb R\,$ as in (\ref{beta}), the processes   $\, \widetilde{\Xi}(\cdot)\,$, $\, \big( \widetilde{X}_{1}(\cdot), \widetilde{X}_{2}(\cdot)\big)\,$ and $\, \big(\widetilde{V}(\cdot), \widetilde{W}(\cdot)\big)\,$ by 
\begin{equation} 
\label{eq: 2DnewX}
\, 
\widetilde{\Xi}(\cdot)\, :=\, 
\widetilde{V}(\cdot) + 2( 1- \beta) 
\widehat{L}^{ \widetilde{Y}}(\cdot) \,, \quad 
\widetilde{X}_{1}(\cdot)\, :=\, \frac{ \, \widetilde{\Xi}(\cdot) + \widetilde{Y}(\cdot) \,}{2}\, , \quad 
\widetilde{X}_{2}(\cdot)\, :=\, \frac{ \, \widetilde{\Xi}(\cdot) - \widetilde{Y}(\cdot) \,}{2}\, , 
\end{equation}
\[
{\mathrm d} \widetilde{V}(\cdot) \, :=\, \big( \rho {\bf 1}_{\{ \widetilde{Y}(\cdot) > 0\}} + \sigma {\bf 1}_{\{ \widetilde{Y}(\cdot) \le 0\}} \big) {\mathrm d} \widetilde{B}_{1}(\cdot) +  \big( \rho {\bf 1}_{\{ \widetilde{Y}(\cdot) \le 0\}} + \sigma {\bf 1}_{\{ \widetilde{Y}(\cdot) > 0\}} \big) {\mathrm d} \widetilde{B}_{2}(\cdot) \, , 
\]
\[
{\mathrm d} \widetilde{W}(\cdot) \, :=\, \big( \rho {\bf 1}_{\{ \widetilde{Y}(\cdot) > 0\}} + \sigma {\bf 1}_{\{ \widetilde{Y}(\cdot) \le 0\}} \big) {\mathrm d} \widetilde{B}_{1}(\cdot) -  \big( \rho {\bf 1}_{\{ \widetilde{Y}(\cdot) \le 0\}} + \sigma {\bf 1}_{\{ \widetilde{Y}(\cdot) > 0\}} \big) {\mathrm d} \widetilde{B}_{2}(\cdot) \, , 
\]
\newpage
\noindent
 Then by (\ref{eq: Y}) and (\ref{eq: 2DnewX}) we obtain 
\begin{equation*} 
\begin{split}
{\text{sgn}} ( \widetilde{Y}( \cdot)) {\mathrm d} W(\cdot)  \, &=\, {\text{sgn}} (\widetilde{Y}( \cdot)) \Big[ \big( \rho {\bf 1}_{\{X_{1}(\cdot) > X_{2}(\cdot)\}} + \sigma {\bf 1}_{\{X_{1}(\cdot) \le X_{2}(\cdot)\}}\big) {\mathrm d} B_{1}(\cdot) \\
& \hspace{5cm} {} - \big( \rho {\bf 1}_{\{X_{1}(\cdot) \le X_{2}(\cdot)\}} + \sigma {\bf 1}_{\{X_{1}(\cdot) > X_{2}(\cdot)\}}\big) {\mathrm d} B_{2}(\cdot) \Big] \\
\, &= \, {\text{sgn}} (\widetilde{Y}(\cdot)) \Big[ \big( \rho {\bf 1}_{\{Y(\cdot) > 0\}} + \sigma {\bf 1}_{\{Y(\cdot) \le 0\}}\big) {\mathrm d} B_{1}(\cdot) - \big( \rho {\bf 1}_{\{Y(\cdot) \le 0\}} + \sigma {\bf 1}_{\{Y(\cdot) > 0\}}\big) {\mathrm d} B_{2}(\cdot) \Big] \,,  \\
{\mathrm d} \widetilde{W}(\cdot) \, &= \, \big( \rho {\bf 1}_{\{ \widetilde{Y}(\cdot) > 0\}} + \sigma {\bf 1}_{\{ \widetilde{Y}(\cdot) \le 0\}} \big) {\mathrm d} \widetilde{B}_{1}(\cdot) -  \big( \rho {\bf 1}_{\{ \widetilde{Y}(\cdot) \le 0\}} + \sigma {\bf 1}_{\{ \widetilde{Y}(\cdot) > 0\}} \big) {\mathrm d} \widetilde{B}_{2}(\cdot)  \, \\
&= {\bf 1}_{\{ \widetilde{Y}(\cdot) > 0\}} \big( \rho\, {\mathrm d} \widetilde{B}_{1}(\cdot) - \sigma {\mathrm d} \widetilde{B}_{2}(\cdot)\big) + {\bf 1}_{\{ \widetilde{Y}(\cdot) \le 0\}} \big( \sigma {\mathrm d} \widetilde{B}_{1}(\cdot) - \rho\, {\mathrm d} \widetilde{B}_{2}(\cdot)\big) \, . 
\end{split}
\end{equation*}
Because of the relationship   
between $\,(B_{1}(\cdot), B_{2}(\cdot))\,$ and $\, (\widetilde{B}_{1}(\cdot), \widetilde{B}_{2}(\cdot))\,$, it can be shown that 
\begin{equation}\label{eq: intertwine}
 {\mathrm d} \widetilde{W}(\cdot) \, =\,  {\text{sgn}} \big( \widetilde{Y}( \cdot)\big)\, {\mathrm d} W(\cdot) \, . 
\end{equation}
In fact, these identities   can be verified formally via the following table: 
 \begin{center}
\begin{tabular}{|c|rr|rcl|}  \hline 
signs of $\, (Y(\cdot), \widetilde{Y}(\cdot)) \,$  & $\, {\mathrm d} \widetilde{B}_{1}(\cdot)\,$ & $\, {\mathrm d} \widetilde{B}_{2}(\cdot) \,$ & $\, {\mathrm d} \widetilde{W}(\cdot) \, $ & $\, = \,$ & $\, {\text{sgn}} ( \widetilde{Y}( \cdot)) {\mathrm d} W(\cdot)\, $\\ \hline  \hline 
$\, (+, + )\,$ &  $\, {\mathrm d} B_{1}(\cdot) \,$ & $\, {\mathrm d}B_{2}(\cdot)\,$ & $\, \rho\,  {\mathrm d} \widetilde{B}_{1}(\cdot) - \sigma \,{\mathrm d} \widetilde{B}_{2}(\cdot)\,$ & $\,  =\,  $ & $\, \rho \,{\mathrm d} B_{1}(\cdot) - \sigma\, {\mathrm d} B_{2}(\cdot)\,$ \\ 
$\,(-, +)\,$  & $\, - {\mathrm d} B_{2}(\cdot)\,$ & $\, - {\mathrm d} B_{1}(\cdot)\, $ & $\,\rho\, {\mathrm d} \widetilde{B}_{1}(\cdot) - \sigma\, {\mathrm d} \widetilde{B}_{2}(\cdot) \,$ & $\,  =\, $ & $\, \sigma \,{\mathrm d} B_{1}(\cdot) - \rho \,{\mathrm d} B_{2}(\cdot) \,$ \\
$\,(+, -)\,$ & $\, {\mathrm d} B_{2}(\cdot) \,$& $\, {\mathrm d} B_{1}(\cdot)\,$ & $\,  \sigma \,{\mathrm d} \widetilde{B}_{1}(\cdot) - \rho \,{\mathrm d} \widetilde{B}_{2}(\cdot) \,$ & $\,  =\, $ & $\,- \rho \,{\mathrm d} B_{1}(\cdot) + \sigma \,{\mathrm d} B_{2}(\cdot) \,$ \\
$\,(-, -)\,$  & $\, - {\mathrm d} B_{1}(\cdot)\,$ & $\, - {\mathrm d} B_{2}(\cdot)\,$& $\,  \sigma\, {\mathrm d} \widetilde{B}_{1}(\cdot) - \rho \,{\mathrm d} \widetilde{B}_{2}(\cdot) \, $ & $\, =\, $ & $\, - \sigma \,{\mathrm d} B_{1}(\cdot) + \rho\, {\mathrm d} B_{2}(\cdot) \,$  \\ \hline 
\end{tabular}
\end{center}	

\medskip
Substituting this relation (\ref{eq: intertwine}) into (\ref{eq: 2Dunfold}) and recalling (\ref{eq: 2DnewX}), we obtain 
\begin{equation} \label{eq: 2DYt}
 {\mathrm d} \big( \widetilde{X}_{1}( t) - \widetilde{X}_{2}( t)\big) \, =\, {\mathrm d} \widetilde{Y}(t) \, =\, {\mathrm d}  \widetilde{W}(t) + 2 (2 \alpha - 1)\, {\mathrm d} \widehat{L}^{\, \widetilde{Y}}( t) \, . 
\end{equation}
Moreover, because of the correspondence between $\,(\widetilde{V}(\cdot), \widetilde{W}(\cdot))\,$ and $\,(V(\cdot), W(\cdot))\,$ and the relation (\ref{eq: 2DnewX}), we obtain 
\begin{equation} \label{eq: sumdiffVW}
\frac{1}{\,2\,} \, {\mathrm d} \big( \widetilde{V}( t) + \widetilde{W}( t) \big) \, =\, \big( \rho {\bf 1}_{\{ \widetilde{Y}( t) > 0\}} + \sigma {\bf 1}_{\{ \widetilde{Y}(  t) \le 0\}} \big)\, {\mathrm d} \widetilde{B}_{1}( t) \, ,  
\end{equation}
\[
\frac{1}{\,2\,} \, {\mathrm d} \big( \widetilde{V}( t) - \widetilde{W}( t) \big) \, =\, \big( \sigma {\bf 1}_{\{ \widetilde{Y}( t) > 0\}} + \rho {\bf 1}_{\{ \widetilde{Y}( t) \le 0\}}\big) \, {\mathrm d} \widetilde{B}_{2}( t) \, .  
\]

 \medskip
 \noindent
Therefore, by calculating the coefficients in front of the local time terms
 and by combining (\ref{eq: 2DnewX}), (\ref{eq: 2DYt}) and (\ref{eq: sumdiffVW}), we can verify that $\, (\widetilde{X}_{1}(\cdot), \widetilde{X}_{2}(\cdot))\, $ satisfies 
\begin{equation} \label{eq: noDskew2}
\begin{split}
{\mathrm d} \widetilde{X}_{1}(t) \, &=\,  \big( \rho {\bf 1}_{\{ \widetilde{X}_{1}(t) > \widetilde{X}_{2}(t)\}} + \sigma {\bf 1}_{\{ \widetilde{X}_{1}(t) \le \widetilde{X}_{2}(t)\}}\big) {\mathrm d} \widetilde{B}_{1}(t)  \, 
+ (2\alpha - \beta) {\mathrm d} \widehat{L}^{\, \widetilde{Y}}(t) \, , \\
{\mathrm d} \widetilde{X}_{2}(t) \, &=\, \big( \rho {\bf 1}_{\{ \widetilde{X}_{1}(t) \le \widetilde{X}_{2}(t)\}} + \sigma {\bf 1}_{\{ \widetilde{X}_{1}(t) > \widetilde{X}_{2}(t)\}}\big) {\mathrm d} \widetilde{B}_{2}(t) 
+ ( 2- 2\alpha - \beta) {\mathrm d} \widehat{L}^{\, \widetilde{Y}}(t) \,  
\end{split}
\end{equation}
that is, (\ref{eq: noDskew}) with the new Brownian motion $\,(\widetilde{B}_{1}(\cdot), \widetilde{B}_{2}(\cdot))\,$. 

\smallskip
By the \textsc{Girsanov}  theorem, we obtain (\ref{eq: 2Dskew3})-(\ref{eq: 2Dskew4}); whereas  the relationship (\ref{eq: 2Dlt}) between the left local time $\,L^{- \widetilde{Y}}(\cdot)\,$ and the right local time $\,L^{ \widetilde{Y}}(\cdot)\,$ allows us now to recover the dynamics of (\ref{eq: 2Dskew1})-(\ref{eq: 2Dskew2}) from those of (\ref{eq: 2D1})-(\ref{eq: 2D2}). \qed

\section{Appendix: Proof of Proposition \ref{Prop5}}
\label{PfProp}

Given a planar Brownian motion $\,(B_{1}(\cdot), B_{2}(\cdot))\,$ on a probability space $\,(\Omega, \mathcal F, \mathbb P) \,$ and real constants $\,\alpha \in (0, 1)\,$, $\,x_{0} \in \R\,$,  we shall construct a process $\,X(\cdot)  :=  \mathfrak q(Y(\cdot))\,$ from the solution $\,Y(\cdot) \,$ of the stochastic differential equation 
\begin{equation} 
\label{eq: strong 1}
Y(\cdot) = \mathfrak p(x_{0}) + \int^{\, \cdot}_{0} \mathfrak s\big(Y(t)\big) \,\mathrm d \big(B_{1}(t) + B_{2}(t)\big)\, ,   
\end{equation}
where $\,\mathfrak p(\cdot)\,$, $\,\mathfrak q(\cdot)\,$ and $\,\mathfrak s(\cdot)\,$ are defined by 
\[
\mathfrak p(x) \, :=\, (1-\alpha)\,  x \, {\bf 1}_{(0, \infty)}(x) + \alpha \, x \, {\bf 1}_{(-\infty, 0]}(x)\, , 
\quad 
\mathfrak q(x) \, :=\, \frac{1}{\, 1-\alpha\, } \, {\bf 1}_{(0, \infty)}(x) + \frac{1}{\, \alpha\, } \, {\bf 1}_{(-\infty, 0)}(x) \, , 
\]
\[
\mathfrak s(x) \, :=\, (1-\alpha)\, {\bf 1}_{(0, \infty)}(x) + \alpha\,  {\bf 1}_{(-\infty, 0]}(x)\, ; \quad x \in \R\, . 
\]
From the work on \textsc{Nakao} (1972) we know that the equation (\ref{eq: strong 1}) has a pathwise unique, strong solution. 

Since $\,\mathfrak q (\mathfrak p(x)) = x\, $, $\,x \in \R\,$, by applying the \textsc{It\^o-Tanaka} formula to the process $\,X(\cdot) = \mathfrak q(Y(\cdot))\,$  we identify the dynamics of $\,X(\cdot)\,$ as those of the skew Brownian motion (\textsc{Harrison \& Shepp} (1981)), namely 
\begin{equation} 
\label{B1_B2}
  X(\cdot) = x_0 +     \big(B_{1}(\cdot) + B_{2}(\cdot)\big) + \frac{\, 2 \alpha - 1\, }{\alpha}   L^{X}(\cdot) \, , 
\end{equation}
driven by the Brownian motion $\,B_{1}(\cdot) + B_{2}(\cdot)\,$.   
We rewrite this equation in the form  
$$
X(\cdot) - x_{0} - \int^{\, \cdot}_{0} \text{sgn}\big(X(t)\big)   \mathrm d U(t) - V(\cdot) \,=\, \frac{  2\alpha - 1 }{\alpha} L^{X}(\cdot)   
\,=\, 2  \big(2\alpha - 1\big) \widehat{L}^{X}(\cdot) 
$$
of (\ref{skewProk}), driven by a new planar Brownian motion $\,( U(\cdot), V(\cdot))\,$ with components
\begin{equation}
\label{U}
U(\cdot) \, :=   \int^{\, \cdot}_{0} {\bf 1}_{\{X(t) > 0\}} \mathrm d B_{1}(t) - \int^{\, \cdot}_{0} {\bf 1}_{\{X(t) \le 0\}} \mathrm d B_{2}(t) \, , 
\end{equation}
\begin{equation}
\label{Q}
V(\cdot) \, :=  \int^{\, \cdot}_{0} {\bf 1}_{\{X(t) \le 0\}} \mathrm d B_{1}(t) + \int^{t}_{0} {\bf 1}_{\{X(t) > 0\}} \mathrm d B_{2}(t)    \, . 
\end{equation}
\noindent
Therefore, the perturbed skew \textsc{Tanaka} equation (\ref{skewProk}) has the   weak solution $\,(X(\cdot), (U(\cdot), V(\cdot)))\,$ just constructed. 

Conversely, suppose we start with an arbitrary weak solution $\,(X(\cdot), (U(\cdot), V(\cdot)))\,$ of the equation (\ref{skewProk}), with  $\,( U(\cdot), V(\cdot))\,$ a  planar Brownian motion. Then we can cast this equation in the form (\ref{B1_B2}) in terms of the  planar Brownian motion $\,( B_1(\cdot), B_2(\cdot))\,$ whose components are given by ``disentangling" in (\ref{U}), (\ref{Q}), namely 
$$
B_1(\cdot) \,  =   \int^{\, \cdot}_{0} {\bf 1}_{\{X(t) > 0\}} \mathrm d U(t) + \int^{\, \cdot}_{0} {\bf 1}_{\{X(t) \le 0\}} \mathrm d V(t) \, , 
$$
$$
B_2(\cdot) \,  =   \int^{\, \cdot}_{0} {\bf 1}_{\{X(t) > 0\}} \mathrm d V(t) - \int^{\, \cdot}_{0} {\bf 1}_{\{X(t) \le 0\}} \mathrm d U(t) \,. 
$$
But this shows that  $\, X(\cdot)\,$ is   skew Brownian motion, so its probability distribution is determined uniquely. 

In other words,  the equation of (\ref{skewProk}) admits a weak solution, and this solution in unique in the sense of the probability distribution. 

 \medskip 
\noindent $\bullet~$ {\it 
Now we shall see that  we have not just uniqueness in distribution, but also   
pathwise uniqueness,  for the equation (\ref{skewProk}) driven by the planar Brownian motion $\,( U(\cdot), V(\cdot))\,$.} The  argument that follows is based on Lemma 1 of \textsc{Le Gall} (1983), and is almost identical to the proof of Theorem 8.1 of [FIKP]  except for the evaluation of the additional local times. Note that \textsc{Le Gall}'s Lemma 1 works for  continuous semimartingales, in general. 

Suppose that there are two solutions $\,X_{1}(\cdot)\,$ and $\,X_{2}(\cdot)\,$ of (\ref{skewProk}), defined on the same probability space as the driving planar Brownian motion $\,( U(\cdot), V(\cdot))\,$. We shall check their difference $\,D(\cdot) \, :=\, X_{1}(\cdot) - X_{2}(\cdot)\,$ satisfies (c.f.  (8.4) in \textsc{Fernholz et al.} (2011)) : 
\begin{equation} \label{eq: Lemma1 LeGall}
\mathbb E \Big[ \int^{T}_{0} \frac{\mathrm d \langle D\rangle(s)}{D(s)} \,{\bf 1}_{\{D(s) > 0 \}} \Big] < \infty \, , \quad 0 < T < \infty \, , 
\end{equation}
where 
\[
\langle D\rangle(\cdot) \, =\, \int^{\cdot}_{0} \big( \text{sgn}(X_{1}(t)) - \text{sgn}(X_{2}(t)) \big)^{2} \mathrm d t \, \le  \, 2 \int^{\cdot}_{0} \big| \text{sgn}(X_{1}(t)) - \text{sgn}(X_{2}(t)) \big|\, \mathrm d t \, .  
\] 
We approximate the  signum function by a sequence $\,\{ f_{k} \}_{k \in \mathbb N} \subset C^{1}(\mathbb R)\,$ which converges to the function $\,f_{\infty}(\cdot) \, =\, \text{sgn}(\cdot)\,$ pointwise and satisfies $\,\lim_{k\to \infty} \lVert f_{k}  \rVert_{TV} = \lVert   f_{\infty} \rVert_{TV}\,$. Now the parametrized process $$\,Z^{(u)}(t) \, :=\, (1-u) X_{1}(t) + u X_{2}(t)\, , \qquad \,0 \le u \le 1\, ,  ~~~\,0 \le t < \infty \,$$ takes the form of 
\[
Z^{(u)}(\cdot) \,= \,x_{0} + \int^{\, \cdot}_{0} \big((1-u) \text{sgn}(X_{1}(t)) + u \, \text{sgn}(X_{2}(t)) \big) \,\mathrm d U(t) 
\]
\[
~~~~~~~~~~~~~~~~~~~~+ \,V(\cdot) + \frac{\, 2 \alpha - 1\, }{\alpha} \big( u L^{X_{1}}(\cdot) + (1-u) L^{X_{2}}(\cdot) \big) \, . 
\]
The local times in the last term do not affect the size of $\,\langle Z^{(u)}\rangle(\cdot)\,$, for which we have the estimate $\, \mathbb{E}  \big(\langle Z^{(u)}\rangle(T) \big)\le 2 T\,$. Proceeding  as in [FIKP] we obtain   for every $\,\delta > 0\,$ the bound 
\[
\mathbb E \Big[ \int^{T}_{0} \frac{\, \lvert f_{k}(X_{1}(s)) - f_{k}(X_{2}(s)) \rvert\, }{X_{1}(s) - X_{2}(s)} {\bf 1}_{\{ X_{1}(s) - X_{2}(s) > \delta\}} \mathrm d t \Big] \, \le\,   c \, \lVert f_{k} \rVert_{TV} \cdot \sup_{a, u} \mathbb E \big( 2 L^{(u)} (T, a) \big) \, , 
\]
 where $\,L^{(u)}(T, a)\,$ is the right local time of the continuous semimartingale $\, Z^{(u)}(\cdot)\,$ accumulated at $\,a \in \R\,$ and $\,c\,$ is a constant chosen independently of $\, k, u, \delta \,$. Letting $\,k \uparrow \infty\,$ and $\,\delta \downarrow 0\,$, we estimate 
\[
\mathbb E \Big[ \int^{T}_{0} \frac{\mathrm d \langle D\rangle(s)}{D(s)} {\bf 1}_{\{D(s) > 0 \}} \Big] \,<\, 2 \,c \, \lVert f_{\infty} \rVert_{TV}  \cdot \sup_{a, u} \mathbb E \big( 2 L^{(u)} (T, a) \big) \, .  
\]
Finally, we estimate $\,\mathbb E(L^{(u)}(T, a)) \,$ using \textsc{Tanaka}'s formula 
\[
\lvert Z^{(u)}(T) - a \rvert = \lvert Z^{(u)}(0) - a \rvert + \int^{T}_{0} \text{sgn} \big(Z^{(u)}(t) - a\big) \mathrm d Z^{(u)}(t) + 2 L^{(u)}(T, a) \, ,  
\]
and a combination of the \textsc{Cauchy-Schwartz} inequality and the \textsc{It\^o}'s isometry:   
\[
\hspace{-6cm} \mathbb E \big(2 L^{(u)}(T, a) \big) \le \mathbb E \lvert Z^{(u)}(T) - Z^{(u)}(0) \rvert + \big\{ \mathbb E ( \langle Z^{(u)}\rangle (T) ) \big\}^{1/2} 
\]
\[
\hspace{4cm} + \frac{\, 2 \alpha - 1\, }{\alpha} \big( u \, \mathbb E(L^{X_{1}}(T) ) + (1-u
)\, \mathbb E( L^{X_{2}}(T)) \big) \]
\[
\hspace{0.5cm} \le 2 \Big[ \big\{ \mathbb E \big( \langle Z^{(u)}\rangle (T) \big) \big\}^{1/2} + \frac{\, 2\alpha -1 \, }{\alpha }
\big( u \, \mathbb E(L^{X_{1}}(T) ) + (1-u
) \,\mathbb E( L^{X_{2}}(T)) \big) \Big] \, . 
\]
The last term $\,\mathbb E(L^{X_{i}}(T))\,$ is evaluated by the same procedure: by \textsc{Tanaka}'s formula 
\[
\frac{1}{\, \alpha\, } L^{X_{i}}(T) \, =\, \lvert X_{i}(T) \rvert - \lvert X_{1}(0) \rvert - \int^{T}_{0} \text{sgn}\big(X_{i}(t)\big)\, \mathrm d V(t) - U(T) \, , 
\]
and hence 
\[
\mathbb E\big( L^{X_{i}}(T)\big) \le 2  \alpha \big\{ \mathbb E ( \langle X_{i} \rangle (T) ) \big\}^{1/2} \le 2^{3/2} \, \alpha \,  T^{1/2} \, , \quad i = 1, 2\, . 
\]
Therefore, we obtain (\ref{eq: Lemma1 LeGall}), and  by Lemma 1 of \textsc{Le Gall} (1983) we verify $\,L^{D}(\cdot)= L^{X_{1} - X_{2}}(\cdot) \equiv 0\,$. 

\medskip 
\noindent $\bullet~$ {\it Final step}: By exchanging the r\^oles of $\,X_{1}(\cdot)\,$ and $\,X_{2}(\cdot)\,$, we obtain $\,L^{-D}(\cdot) = L^{X_{2} - X_{1}}(\cdot) \equiv 0\,$ as well as $\, \widehat{L}^{D}(\cdot) \equiv 0\,$. 
Furthermore, by Corollary 2.6 of \textsc{Ouknine \& Rutkowski} (1995), we obtain   
\[
\widehat{L}^{X_{1} \vee X_{2}}(t) \, =\, \int^{t}_{0} 
{\bf 1}_{\{X_{2}(s) \le 0\}}\, {\mathrm d} \widehat{L}^{X_{1}}(s) + \int^{t}_{0} {\bf 1}_{\{X_{1}(s) < 0\}} \,{\mathrm d} \widehat{L}^{X_{2}}(s) \, ; \quad 0 \le t < \infty \, . 
\]
Combining these results with \textsc{Tanaka}'s formula, we obtain the dynamics of $\,M(\cdot) \, :=\, X_{1}(\cdot) \vee X_{2}(\cdot)\,$: 
\[
{\mathrm d} M(t) \, =\, {\bf 1}_{\{X_{1} (t) \ge X_{2}(t)\}}{\mathrm d} X_{1}(t) + {\bf 1}_{\{X_{1}(t) < X_{2}(t)} {\mathrm d} X_{2} (t) + {\mathrm d} L^{X_{1} - X_{2}}(t) ~~~~~~~~~~~~~~~~~~ 
\]
\[
~~\, =\, {\bf 1}_{\{X_{1} (t) \ge X_{2}(t)\}} \Big( \text{sgn} (X_{1}(t)) {\mathrm d} U(t) +  {\mathrm d} V(t) + 2(2\alpha - 1) {\mathrm d} \widehat{L}^{X_{1}}(t) \Big) \]
\[
\hspace{3cm} {}+ {\bf 1}_{\{X_{1}(t) < X_{2}(t)}\Big( \text{sgn} (X_{2}(t)) {\mathrm d} U(t) +  {\mathrm d} V(t) + 2(2\alpha - 1) {\mathrm d} \widehat{L}^{X_{2}}(t) \Big) \, 
\]
\[
\, =\,  \text{sgn} (M(t)) {\mathrm d} U(t) +  {\mathrm d} V(t) + 2(2\alpha - 1)  {\mathrm d} \widehat{L}^{M}(t)  \, ; \quad 0 \le t < \infty \, . 
\]

\medskip
In other words, each of the continuous semimartingales $\,X_{1}(\cdot)\,$, $\,X_{2}(\cdot)\,$ and $\,M(\cdot) = X_{1}(\cdot) \vee X_{2}(\cdot)\,$ satisfies the equation (\ref{skewProk}); but uniqueness in the sense of the probability distribution holds for this equation, so all three processes have the same distribution. Since $\, M(\cdot)\ge X_i (\cdot)\,$, this forces  $\, M(\cdot)= X_i (\cdot)\,$, $\,i=1, 2\,$, thus    pathwise uniqueness. By the theory of \textsc{Yamada} and \textsc{Watanabe} (e.g., subsection 5.3.D in \textsc{Karatzas \& Shreve} (1991)), the solution to (\ref{skewProk}) is   therefore strong.  
The proof of Proposition \ref{Prop5} is complete. \qed



\section*{Bibliography}


\noindent \textsc{Blei, S.}  (2012) On symmetric and skew Bessel processes. {\it Stochastic Processes and Their Applications} {\bf 122}, no.~9, 3262--3287. 

 \medskip 
  \noindent
  \textsc{Chaleyat-Maurel, M. \& Yor, M.} (1978) 
 Les filtrations de $\,|X|\,$ et $\,X^+\,$, lorsque $\,X\,$ est une semi-martingale continue.   In   ``Temps Locaux".  {\it  Ast\'erisque} {\bf 52-53},  193-196.

\medskip 
\noindent \textsc{Cherny, A. S.} (2000) On the strong and weak solutions of stochastic differential equations governing Bessel processes, {\it Stochastics \& Stochastics Reports} {\bf 70}, no.~3-4, 213--219. 

\medskip 

\noindent \textsc{Dubins, L.E., \'Emery, M. \& Yor, M.} (1993) On the L\'evy transformation of Brownian motions and continuous martingales. In ``S\'eminaire de   Probabilit\'es XXVII ". {\it Lecture Notes in Mathematics} {\bf 1577}, 122-132. Springer Verlag, New York.

\medskip 

\nocite{FIK12}
\noindent \textsc{Fernholz, E.R., Ichiba, T. \& Karatzas, I.} (2013) [FIK] {Two {B}rownian particles with rank-based characteristics and skew-elastic collisions.} {\it Stochastic Processes and Their Applications } {\bf 123}, 2999-3026.

 \medskip 

\nocite{FIKP}
\noindent \textsc{Fernholz, E.R., Ichiba, T., Karatzas, I. \& Prokaj, V.} (2013) [FIKP] {A planar diffusion with rank-based characteristics, and perturbed Tanaka equations.} {\it Probability Theory and Related Fields} {\bf 156}, 343-374.

\medskip
\noindent\textsc{Harrison, J.M. \& Shepp, L.A.} (1981). On skew Brownian motion. {\it Annals of Probability} {\bf  9}, 309--313.

 \medskip
 \nocite{IKP}
\noindent \textsc{Ichiba, T., Karatzas, I. \& Prokaj, V.} (2013)  {Diffusions with rank-based characteristics and values in the nonnegative orthant.} {\it Bernoulli} {\bf 19}, 2455-2493.

 \medskip
 
 \noindent    \textsc{Karatzas, I. \& Shreve, S.E.}  (1991)
  {\it Brownian Motion and Stochastic Calculus.} Second Edition, Springer Verlag, New York.
  
   \medskip
\noindent \textsc{Le Gall, J-F.} (1983) Applications du temps local aux \'equations diff\'erentielles stochastiques unidimensionnelles. In ``S\'eminaire de   Probabilit\'es XVII". {\it Lecture Notes in Mathematics} {\bf 986}, 15-31. Springer Verlag, New York.

  \medskip
\noindent \textsc{Nakao, S.}  (1972)  On pathwise uniqueness of solutions of stochastic differential equations. {\it Osaka Journal of Mathematics} {\bf  9}, 513-518.

\medskip 
\noindent \textsc{Ocone, D.} (1993) A symmetry characterization of conditionally independent increment martingales. {\it Proceedings of the San Felice Workshop on Stochastic Analysis} (D.$\,$Nualart \& M.$\,$Sanz, editors), 147-167. Birkh\"auser-Verlag, Basel and Boston.

\medskip 
\noindent \textsc{Ouknine, Y. \& Rutkowski, M.} (1995) Local times of functions of continuous
semimartingales. {\it Stochastic Analysis and Applications} {\bf 13}, 211-231. 
 
   \medskip
\noindent \textsc{Prokaj, V.}  (2009)  Unfolding the Skorokhod reflection of a semimartingale. {\it Statistics and  Probability Letters} {\bf 79}, 534-536.

  \medskip
\noindent \textsc{Prokaj, V.}  (2013)  The solution of the perturbed Tanaka equation is pathwise unique. {\it Annals of Probability} {\bf 41}, 2376-2400.

   \medskip
\noindent \textsc{Revuz, D. and Yor, M.} (1999) {\it Continuous Martingales and Brownian Motion}. Third edition, Springer Verlag, New York.

 \medskip
\noindent \textsc{Vostrikova, L. \& Yor, M.} (2000) Some invariance properties of the laws of Ocone  martingales. In ``S\'eminaire de   Probabilit\'es XXXIV ". {\it Lecture Notes in Mathematics} {\bf 1729}, 417-431. Springer Verlag, New York.

\end{document}